\documentclass[10pt,a4paper]{article}
\usepackage{amsmath,amssymb,amsfonts,amsthm,graphics,graphicx,psfrag}
\newcommand{\eps}{\varepsilon}
\newcommand{\R}{\mathbb{R}}

\renewcommand{\a }{\alpha }
\renewcommand{\le }{\leqslant }
\renewcommand{\ge }{\geqslant }
\renewcommand{\b }{\beta }

\newcommand{\G}{\Gamma}

\newtheorem{theorem}{Theorem}[section]
\newtheorem{lemma}[theorem]{Lemma}
\newtheorem{claim}[theorem]{Claim}

\newtheorem{proposition}[theorem]{Proposition}

\newtheorem{corollary}[theorem]{Corollary}
\renewenvironment{proof}{\noindent{\textbf{Proof\quad}}}{$\hfill\square$\vspace{0.2
cm}\\}
\newenvironment{proofl}{\noindent{\textbf{Proof of Theorem
\ref{stabl}.\quad}}}{$\hfill\square$\vspace{0.2 cm}\\}

\numberwithin{equation}{section}

\begin{document}
{\title{ \bf{Stable determination of the surface\\
\bf impedance of an obstacle by\\
 \bf  far field measurements}
\thanks{Work supported in part by MIUR, grant n.
2004011204.} }}
\author{E. Sincich \thanks{S.I.S.S.A., Trieste, Italy, \tt {
sincich@sissa.it}}}

\date{}

\maketitle

\begin{abstract}We deal with the inverse scattering problem of
determining the surface impedance of a partially coated obstacle. We
prove a
stability estimate of logarithmic type for the impedance term by the
far
field measurements.
\end{abstract}



\section{Introduction}
We consider the scattering of an acoustic incident time-harmonic plane
wave, at a given wave number $k>0$ and at a given incident direction
$\omega\in \mathbb{S}^2$, by an obstacle $D\subset
\mathbb{R}^3$ partially coated by a material with surface impedance
$\lambda$. Such a problem is modeled by the following mixed boundary
value problem for the Helmholtz equation
\begin{eqnarray}\label{Sc}
\left\{
\begin{array}
{lll}
\Delta u +k^2u=0,\ \ \  &\mbox{in $\R^3 \setminus
\overline{D}$},
\\
u=0,\ \ \ &\mbox{on $\G_D$},\ \ \
\\
\dfrac{\partial u}{\partial \nu} + i\lambda(x)u=0,\ \ \
&\mbox{on $\G_I$},
\end{array}
\right.
\end{eqnarray}
where $u=u^s+\exp{(ikx\cdot \omega)}$ is the total field,
that is given as the sum of the scattered wave $u^s$ and the incident
plane waves $\exp{(ikx\cdot \omega)}$ and where $\G_I,\ \G_D$ are
two open and connected portions of the boundary $\partial D$ such that
$\partial
D=\overline {\G_I\cup\G_D}$.

Moreover, the scattered field $u^s$ is required to satisfy the
so-called \emph{Sommerfeld
radiation condition}
\begin{eqnarray}\label{Sc3}
\ \ \ & \lim_{_{r\rightarrow\infty}}r\left(\dfrac{\partial
u^{s}}{\partial
r}-iku^s \right)=0,\   \ & r=\|x\|.
\end{eqnarray}

It is well-known, that the scattered
field $u^s$ has the following  asymptotic behavior
\begin{eqnarray}
u^s(x)=\frac{\exp{(ikr)}}{r}\left\{u_{\infty}(\hat{x})+O\left(\frac{1}{r}\right)\right\}\
,
\end{eqnarray}
as $r$ tends to $\infty$, uniformly with respect to
$\hat{x}=\frac{x}{\|x\|}$ and where
$u_{\infty}$ is the so-called far field pattern of the scattered wave
(see for instance \cite{ck}).

The \emph{inverse scattering problem} that we examine here consists in
the determination of the
surface impedance $\lambda(x)$ by the knowledge of the far field
pattern, provided some suitable
a priori assumptions on the impedance are made.

Such a problem, in two dimensions, has been recently studied
by F. Cakoni and D. Colton in \cite{cc}. The authors have provided a
variational method for the determination of the essential supremum of
the
surface impedance when the far field data are available.

In this paper, we shall deal with the stability issue, namely
we will prove a stability estimate of logarithmic type for the surface
impedance by the far field measurements.

Let us point out that a stability result for this type of problem has
been proved in \cite{francese} by C. Labreuche under the assumption of
an analytic boundary. The new feature of the present paper consists in
a
reduced assumption on the regularity of the boundary, namely we shall
assume that $\G_I$ is a $C^{1,1}$ portion of $\partial D$. Thus it
turns
out that the argument of analytic extension used in \cite{francese}
cannot be applied.

The stable recovering of the surface impedance needs some a priori mild
assumptions on the impedance itself.
The additional a priori information that we require on the unknown
surface impedance $\lambda$, is an a priori bound on its Lipschitz
continuity, that is we assume that for a given positive constant
$\Lambda$, the following holds
\begin{eqnarray}\label{regintr}
\|\lambda\|_{C^{0,1}(\G_I)}\le \Lambda .
\end{eqnarray}
Moreover, we prescribe the following uniform lower bound
\begin{eqnarray}\label{dalbassointr}
\ \ \ \ \ \ \ \ \ \ \ \ \  \ \ \ \ \ \lambda(x)\ge\lambda_0 ,\ \
\mbox{for every}\ x\in\G_I,
\end{eqnarray}
where $\lambda_0$ is a given positive constant.

In order to treat the inverse scattering problem we first need to
analyze the \emph{direct} one. In Section 3, indeed, following the
arguments of potential theory developed in \cite{ccm}, we observe that
the direct
scattering problem is well posed (see Lemma \ref{benposto}). The proof
relies on the fact that the mixed boundary value problem \eqref{Sc} can
be reformulated as a system of boundary integral equations. Moreover,
we prove, (see Theorem \ref{regolarita}), that the solution and its
first order derivatives are H\"{older} continuous in a neighborhood of
the
portion $\G_I$, where the impedance takes place. The proof is based on
the Moser's iteration technique. Finally in Corollary \ref{lb}, we
obtain a uniform lower bound for the total field $u$ on sets away from
the
obstacle.

In Section 4, we deal with the inverse scattering problem. The
underlying
ideas and the main tools that lead to the stability result can be
outlined as follows.
\begin{description}
\item[i)]As first step we evaluate how much the error on the far field
can
affect the values of the field near the scatterer;
\item[ii)]in the second step we are concerned with a stability estimate
of the
field at the boundary in terms of the near field;
\item[iii)]finally, as last step, we obtain a stability result for the
impedance $\lambda$ by the estimate of the field at the boundary.
\end{description}
Let us start the analysis of Section 4 illustrating the arguments
introduced in the step iii) of the list above.

By the impedance condition in \eqref{Sc} we can formally compute
$\lambda$ as
\begin{eqnarray}\label{formalmente}
\lambda(x)=\frac{i}{u(x)}\dfrac{\partial u(x)}{\partial \nu(x)}\ .
\end{eqnarray}
Since $u$ may vanish in some points of $\G_I$, it follows that
the quotient in \eqref{formalmente} may be undetermined. In this 
respect, we found it
necessary to evaluate the local vanishing rate of the solution on the
boundary.
To establish such a control we shall make use of quantitative estimates
of unique continuation. We first obtain, in  Lemma \ref{dib}, a
\emph{volume doubling inequality} at the boundary, namely
\begin{eqnarray}\label{intrat}
\int_{\G_{I,2\rho}(x_0)}|u|^2\le
\emph{const.}\int_{\G_{I,\rho}(x_0)}|u|^2\ ,
\end{eqnarray}
where $\G_{I,\rho}(x_0)$ and $\G_{I,2\rho}(x_0)$ are the portions of
the balls centered at the boundary point $x_0$ of radius $\rho$ and
$2\rho$ respectively, contained in ${\R}^3\setminus \overline{D}$, (see
\eqref{neigh} for a precise definition).

In order to obtain the formula in \eqref{intrat}, we have
adapted the arguments developed in \cite{ae} for the more general
setting of complex valued solutions which is required by the boundary
value problem \eqref{Sc}.

 A further difficulty in dealing with such arguments is due to the fact
that the techniques used in \cite{ae} apply
to an homogeneous Neumann boundary condition. We overcome such a
difficulty by performing a suitable change of the independent variable,
(see
Proposition \ref{psi2}), that fits our problem under the assumptions
required in \cite{ae}. Moreover, well-known stability estimates for the
Cauchy problem \cite{tr}, allow us to reformulate the \emph{volume} 
doubling inequality at the boundary deriving in Theorem \ref{di}
a new one on the boundary, that is a \emph{surface doubling inequality}
\begin{eqnarray}\label{inton}
\int_{\Delta_{I,2\rho}(x_0)}|u|^2\le
\emph{const.}\int_{\Delta_{I,\rho}(x_0)}|u|^2\ ,
\end{eqnarray}
where $\Delta_{I,\rho}(x_0)$ and $\Delta_{I,2\rho}(x_0)$ are the
portions of the boundary of $\G_{I,\rho}(x_0)$ and $\G_{I,2\rho}(x_0)$
respectively, which have non empty intersection with  $\partial D$,
(see \eqref{neighb} for a precise
definition).

The surface doubling inequality allows us to apply the theory of
\emph{Muckenhoupt weights}
\cite{cf} which, in particular, implies the existence of some exponent
$p>1$ such that $|u|^{-\frac{2}{p-1}}$ is integrable on an inner
portion
of $\G_I$, see Corollary \ref{ap}. This integrability property, as well
as the H\"{o}lder continuity of the normal derivative, justifies the
computation made in \eqref{formalmente} in the $L^{\frac{2}{p-1}}$
sense.

Let us carry over our analysis by discussing the evaluation
introduced in the step i).
Such an evaluation, introduced by V. Isakov \cite{is,is2}, and then
developed by I. Bushuyev \cite{b}, concerns a stability estimate for
the \emph{near field} in terms of the measurements of the \emph{far
field} (see Lemma
\ref{farnear}). It means that if $u_1$ and $u_2$ are two acoustic
fields corresponding to impedances $\lambda_1$ and $\lambda_2$ such
that
their scattering amplitudes, $u_{1,\infty}$ and $u_{2,\infty}$
respectively, are close
\begin{eqnarray}\label{intrerr}
\|u_{1,\infty}-u_{2,\infty}\|_{L^2(\partial B_1(0))}\le \varepsilon,
\end{eqnarray}
then $u_1$ and $u_2$ satisfy
\begin{eqnarray}\label{intrstab1}
\|u_1-u_2\|_{L^2(B_{R_1+1}(0)\setminus B_{R_1}(0))}\le \emph{const.}
{\varepsilon}^{\alpha(\eps)},
\end{eqnarray}
where $R_1>0$ is a suitable radius such that $B_{R_1}(0)\supset
\overline{D}$ and $\alpha(\eps)$ is the function introduced in
\eqref{lam}.

As last step of this treatment we provide the stability estimate
introduced in ii).
 The proof is based on arguments of quantitative unique continuation,
as the \emph{three
spheres inequality} and leads to the following estimate
\begin{eqnarray}\label{intrstab2}
\|u_1-u_2\|_{C^1(\G^{\rho}_I)}\le
\emph{const.}|\log{(\|u_1-u_2\|^{-1}_{L^2(B_{R_1+1}(0)\setminus
B_{R_1}(0))})|^{-\theta}},
\end{eqnarray}
 where $\theta>0$ and where $\G^{\rho}_I$ is a given inner portion of
$\G_I$ (see \eqref{innerportion} for a precise definition).

By combining the stability estimates listed in i) and ii), we obtain
a stability result for the total field at the boundary in terms of the
measurements of the far field, (see Theorem \ref{stabpc}).

Finally, as a consequence of Theorem \ref{stabpc} and Corollary
\ref{ap}, let us formulate the main result of the present paper, that
consists
in a stability estimate of
the surface impedance by the far field measurements, (see Theorem
\ref{stabl}). Assuming that \eqref{intrerr} holds, we have shown that
the impedances $\lambda_1, \lambda_2$ agree up to an error
\begin{eqnarray}\label{intrris}
\big|\log\big({\varepsilon}^{-\alpha(\eps)}\big)\big|^{-\theta}.
\end{eqnarray}
Moreover, let us observe that the rate of stability in \eqref{intrris}
is
intermediate between a \emph{log} and a \emph{loglog} rate of
stability.

\section{Main assumptions and results}
\subsection{Main hypothesis and notations}
\emph{\bf Assumptions on the domain.}\\
We shall assume through that $D$ is a bounded domain in $\R^3$, such
that $\mbox{diam}D\!\!\le\!\!d$, with Lipschitz boundary $\partial D$
with
constants $r_0, M$.
More precisely,
for every $x_0 \in \partial D$, exists a rigid transformation of
coordinates under which,
\begin{eqnarray}\label{gamma}
D \cap B_{r_0}(x_0)=\{(x',x_3): x_3>\gamma(x')\}\ ,
\end{eqnarray}
where $x \in \R^3,\ x=(x',x_3)$, with $x'  \in \R^2 , \ x_3 \in
\mathbb{R}$ and
$$\gamma :B'_{r_0}(x_0)\subset \R^2 \rightarrow \mathbb{R}$$
satisfying $\gamma(0)=0$ and $$\|\gamma\|_{C^{0,1}( B'_{r_0}(x_0))}\le
Mr_0,$$ where we denote by
$$\|\gamma\|_{C^{0,1}( B'_{r_0}(x_0))}=\|\gamma\|_{L^{\infty}(
B'_{r_0}(x_0))} +\ r_0\!\!\!\!\!\!\!\!\!\sup_{\substack {x,y  \in
B'_{r_0}(z_0)\\x\ne y }}
\frac{|\gamma (x)-\gamma (y)|}{|x-y|} _{\ \ \  } $$
and $B'_{r_0}(x_0)$ denotes a ball in $\R^2$.
Moreover, we assume that the portion of the boundary $\Gamma_I$ is
contained into a surface $S_I$, which is $C^{1,1}$ smooth with
constants
$r_0, M$.

More precisely, for any $x_0 \in S_I$, we have that up to a rigid
change of coordinates,
\begin{eqnarray}\label{sI}
  S_I \cap B_{r_0}(x_0)=\{(x',x_3): x_3=\varphi_I(x')\},
\end{eqnarray}
where
\begin{eqnarray}\label{sI1}
\varphi_I:B'_{r_0}(z_0)\subset \R^2 \rightarrow \mathbb{R}
\end{eqnarray}
 is a $C^{1,1}$ function satisfying
\begin{eqnarray}\label{sI2}
\varphi_I(0)=|\nabla \varphi_I(0)|=0
\end{eqnarray}
 and
\begin{eqnarray}\label{sI3}
\|\varphi_I \|_{C^{1,1}(B'_{r_0}(z_0))}\le Mr_0 ,
\end{eqnarray}
where we denote
\begin{eqnarray}\label{sI4}
\|\varphi_I\|_{C^{1, 1}( B'_{r_0}(z_0))}&=&\|\varphi_I\|_{L^{\infty}(
B'_{r_0}(z_0))}+r_0\|\nabla \varphi_I\|_{L^{\infty}( B'_{r_0}(z_0))}+\\
& & +\ {r_0}^{2}\!\!\!\!\!\!\!\!\!\sup_{\substack {x,y  \in
B'_{r_0}(z_0)\\x\ne y }}\frac{|\nabla \varphi_I (x)-\nabla \varphi_I
(y)|}{|x-y|}_{\ \ \ .}
\end{eqnarray}

In particular it follows that, if
$$x_0 \in \Gamma_I\ \ \mbox{and} \ \ \mbox{dist}(x_0,\G_D)>r_0\ ,$$
then
\begin{eqnarray}\label{sI5}
D \cap B_{r_0}(x_0)=\{(x',x_3) \in  B_{r_0}(x_0): x_3>\varphi_I(x')\}\
,
\end{eqnarray}
where $\varphi_I$ is the Lipschitz function whose graph locally
represents $\partial D$. Moreover, since $D \cap B_{r_0}(x_0) \cap \G_D
=\emptyset$, $\varphi_I$ must also be the $C^{1,1}$ function whose
graph
locally represents $S_I$.\\

For a sake of simplicity we shall assume that $0\in D$.

Fixed $R>d$, $\rho\in(0,r_0)$ and $x_0\in\G_I$, let us define the
following sets

\begin{eqnarray}
&&D^{+}=\R^3\setminus\overline{D},\\
&&D_R^{+}=B_R(0)\cap D^{+}, \\
&&D_{R,\rho}^{+}=\{x\in \overline{D_R^{+}}\
:\mbox{dist}(x,\G_D)>\rho\},\\
&&\G_{I}^{\rho}=\partial D_{R,\rho}^{+} \cap \G_I ,\label{innerportion}
\\
&&\G_{I,\rho}(x_0)= B_{\rho}(x_0)\setminus \overline{D}\label{neigh}
,\\
&&\Delta_{I,\rho}(x_0)=\overline{\G_{I,\rho}(x_0)}\cap \partial
D.\label{neighb}
\end{eqnarray}

\emph{\bf A priori information on the impedance term.}\\
We assume that the impedance coefficient $\lambda$ belongs to
$C^{0,1}(\G_I,\R)$ and is such that
\begin{eqnarray}\label{limdalbasso}
\lambda(x)\ge\lambda_0>0
\end{eqnarray}
 for every $x\in \G_I$. Moreover we assume that, for a given constant
$\Lambda>0$, we
have that
\begin{eqnarray}\label{L}
\|\lambda\|_{C^{0,1}(\G_I)}\le \Lambda .
\end{eqnarray}

From now on we shall refer to the \emph{a priori data} as to the
following set of quantities:
$d, r_0, M,\lambda_0, \Lambda, k,\omega$.\\
In the sequel we shall denote with $\eta(t)$ a positive increasing
function defined on $(0, +\infty)$, that satisfies
\begin{eqnarray}\label{eta}
&&\eta(t)\le C(\log (t^{-\alpha(t)}))^{-\vartheta},\ \ \ \mbox{for
every}\ \ 0<t<1 \ ,
\end{eqnarray}
where
\begin{eqnarray}\label{lam}
\a(t)=\frac{1}{1+\log(\log({t}^{-1})+e)}\ ,
\end{eqnarray}
 and $C>0,\theta>0$ are constants depending on the \emph{a priori
data} only.

\subsection{The main result}
\begin{theorem}[Stability for $\lambda$]\label{stabl}
Let $u_i, \ i=1,2$, be the weak solutions to the problem \eqref{Sc}
with $\lambda=\lambda_i$ respectively and let $u_{i,\infty}$ be their
respectively far field patterns.
There exist $\delta>0, \eps_0>0$ constants only depending on the
\emph{a priori data}, such that, if for some
 $\varepsilon,\ 0<\varepsilon<\eps_0$, we have
\begin{eqnarray}\label{farf}
\|u_{1,\infty}-u_{2,\infty}\|_{L^2(\partial B_1(0))}\le \varepsilon,
\end{eqnarray}
then
\begin{eqnarray}\label{aim}
\|\lambda_1-\lambda_2\|_{L^{\infty}(\G^{r_0}_I)}\le \eta(\eps),
\end{eqnarray}
where $\eta$ is given by \eqref{eta}.
\end{theorem}

\section{The direct scattering problem}

Let us introduce the following space
$$H^1_{\mbox{loc}}(D^+)=\{v\in \mbox{\itshape D}^*(D^+):\
v|_{D_{R}^+}\in H^1(D_{R}^+),\ \mbox{for every}\ R>0\ \mbox{s.t.}\
\overline{D}\subset B_{R}(0)\}$$
where $\mbox{\itshape D}^*(D^+)$ is the space of distribution on $D^+$.

A weak solution to the problem \eqref{Sc} is a function
$u=\exp{(ik\omega\cdot x)}+u^s$, where $u^s\in
H^1_{\mbox{loc}}(D^+)$ is a weak solution to the problem
\begin{eqnarray}\label{radiating}
\left\{
\begin{array}
{lll}
\Delta u^s +k^2u^s=0,& \mbox{in $D^+$},
\\
u^s=-\exp{(ik\omega\cdot x)},  & \mbox{on $\G_D$},
\\
\dfrac{\partial u^s}{\partial \nu} +
i\lambda(x)u^s=-\dfrac{\partial}{\partial \nu}\exp{(ik\omega\cdot x)} -
i\lambda(x)\exp{(ik\omega\cdot
x)}, & \mbox{on $\G_I$},
\\
\lim_{r\rightarrow\infty}r\left(\dfrac{\partial u^{s}}{\partial
r}(r\hat{x})-iku^s(r\hat{x}) \right)=0, & \mbox{uniformly in}\
\hat{x}.
\end{array}
\right.
\end{eqnarray}

Let us recall that a weak solution of \eqref{radiating} is a function
$u^s\in H^1_{\mbox{loc}}(D^+)$, with $u^s|_{\G_D}=-\exp{(ik\omega\cdot
x)}$ in the trace sense, such that, for all test functions $\eta \in
H^1(D^+)$ with compact support in $\R^3$ and $\eta|_{\G_D}=0$, the
following holds
\begin{eqnarray}\label{weak}
\int_{D^+}\nabla u^s\cdot\nabla\overline{\eta}-k^2\int_{D^+}
u^s\overline{\eta}&=&\!\!\!\int_{\G_I}\left(\dfrac{\partial}{\partial
\nu}\exp{(ik\omega\cdot x)} + i\lambda(x)\exp{(ik\omega\cdot
x)}\right)\overline{\eta} +\nonumber\\
&&\!\!\!+\int_{\G_I}ik\lambda u^s\overline{\eta}\ .
\end{eqnarray}
Furthermore, $u^s$ satisfies the asymptotic condition \eqref{Sc3}.

\begin{lemma}[\textbf{Well-posedness}]\label{benposto}
The problem \eqref{radiating} has one and only one weak solution $u^s$.
Moreover, for every $R>d$, there exists a constant $C_R>0$ depending
on the \emph{a priori data} and on $R$ only, such that the following
holds
\begin{eqnarray}\label{dipcont}
\|u^s\|_{H^1(D_R^+)}\le C_R\ .
\end{eqnarray}
\end{lemma}
\begin{proof}
For the proof we refer to \cite[Theorem 2.5]{ccm}, in which the
authors, among various results, show that the exterior mixed boundary
value
problem
\eqref{radiating} can be reformulated as a $2\times 2$ system of
boundary
integral equations.
In \cite{ccm}, Theorem 2.5 has been proved in two dimensions for a
constant $\lambda$, however it can be verified that the same techniques
can be carried over in
three
dimensions
 and with
$\lambda=\lambda(x) \in C^{0,1}(\G_I)$.
\end{proof}
\begin{theorem}[\textbf{$C^{1,\a}$ regularity at the
boundary}]\label{regolarita}

Let $u$ be the weak solution to \eqref{Sc}, then there exists a
constant $\a,\ 0<\a<1$, such that for every $R>d$ and $\rho\in
(0,r_0)$, $u\in C^{1,\a}(D_{R,\rho}^+)$. Moreover, there exists a
constant $C_{R,\rho}>0$ depending on the \emph{a priori data}, on
$R$ and on $\rho$ only, such that
\begin{eqnarray}\label{ciunoalfa}
\|u\|_{C^{1,\a}(D_{R,\rho}^+)}\le C_{R,\rho}\ .
\end{eqnarray}
\end{theorem}
\begin{proof}
From the weak formulation \eqref{weak}, it follows that the total field
$u$ satisfies
 \begin{eqnarray*}
\int_{\G_{I,\frac{r_0}{2}}(x_0)}\!\!\!\!\!\!\!\nabla u\cdot \nabla
\bar{\eta}-k^2\int_{\G_{I,\frac{r_0}{2}}(x_0)}\!\!\!\!\!\!u
\bar{\eta}=-i\int_{\Delta_{I,\frac{r_0}{2}}(x_0)}\!\!\!\!\!\!\!\lambda(x)u\bar{\eta}\ ,
\end{eqnarray*}
where $x_0\in \G_{I}$ and $\eta$ is any test function such that
$\mbox{supp}{\eta}\subset {\overline{\G}_{I,\frac{r_0}{2}}(x_0)}$.\\By
\eqref{L} we have that
\begin{eqnarray}
\bigg|\int_{\G_{I,\frac{r_0}{2}}(x_0)}\!\!\!\!\!\!\!\nabla u\cdot
\nabla \bar{\eta}\bigg |\le
k^2\int_{\G_{I,\frac{r_0}{2}}(x_0)}\!\!\!\!\!\!|u \bar{\eta}|
+
\Lambda\int_{\Delta_{I,\frac{r_0}{2}}(x_0)}\!\!\!\!\!\!\!|u\bar{\eta}|
\end{eqnarray}
and by a trace inequality (see \cite[p.114]{Adams}) it follows that
\begin{eqnarray}
\bigg|\int_{\G_{I,\frac{r_0}{2}}(x_0)}\!\!\!\!\!\!\!\nabla u\cdot
\nabla \bar{\eta}\bigg |\le
k^2\int_{\G_{I,\frac{r_0}{2}}(x_0)}\!\!\!\!\!\!|u \bar{\eta}|
+
C\Lambda\int_{\G_{I,\frac{r_0}{2}}(x_0)}\!\!\!\!\!\!\!|\nabla
(u\bar{\eta})|\ ,
\end{eqnarray}
where $C>0$ is a constant depending on the \emph{a priori data} only.

By the standard iteration techniques due to Moser (see for instance
\cite{gt}), we obtain the following local bound for $u$
\begin{eqnarray}\label{lc}
\|u\|_{L^{\infty}(\G_{I,\frac{r_0}{4}}(x_0))}\le
C\|u\|_{H^1\big(\G_{I,\frac{r_0}{2}}(x_0)\big)}\ ,
\end{eqnarray}
where $C>0$ is a constant depending on the \emph{a priori data} only.

Let us denote by $u_1$ and $u_2$ the real and the imaginary part of $u$ 
respectively. Thus by the elliptic
 equations in weak form satisfied by $u_1$ and $u_2$, it follows that
\begin{eqnarray}
\int_{\G_{I,\frac{r_0}{2}}(x_0)}\!\!\!\!\!\!\!\nabla u_1\cdot \nabla
{\eta}-k^2\int_{\G_{I,\frac{r_0}{2}}(x_0)}\!\!\!\!\!\!u_1
{\eta}=\int_{\Delta_{I,\frac{r_0}{2}}(x_0)}\!\!\!\!\!\!\!\lambda(x)u_2{\eta}\ 
,\label{wk1}\\
\int_{\G_{I,\frac{r_0}{2}}(x_0)}\!\!\!\!\!\!\!\nabla u_2\cdot \nabla
{\eta}-k^2\int_{\G_{I,\frac{r_0}{2}}(x_0)}\!\!\!\!\!\!u_2
{\eta}=-\int_{\Delta_{I,\frac{r_0}{2}}(x_0)}\!\!\!\!\!\!\!\lambda(x)u_1{\eta}\ 
,\label{wk2}
\end{eqnarray}
where $\eta$ is any real valued test function such that
$\mbox{supp}{\eta}\subset {\overline{\G}_{I,\frac{r_0}{2}}(x_0)}$.

By applying again the Moser method to the weak formulations \eqref{wk1} 
and \eqref{wk2}, we obtain the
following bounds of the H\"{o}lder continuity of $u_1$ and $u_2$, 
namely
\begin{eqnarray}
\|u_1\|_{C^{0,\a}(\G_{I,\frac{r_0}{8}}(x_0))}\le
C(\|u_1\|_{L^{\infty}(\G_{I,\frac{r_0}{4}}(x_0))}+\|u_2\|_{L^{\infty}(\G_{I,\frac{r_0}{4}}(x_0))} 
)\ ,\label{cunoalfa1}\\
\|u_2\|_{C^{0,\a}(\G_{I,\frac{r_0}{8}}(x_0))}\le
C(\|u_2\|_{L^{\infty}(\G_{I,\frac{r_0}{4}}(x_0))}+\|u_1\|_{L^{\infty}(\G_{I,\frac{r_0}{4}}(x_0))} 
)\ ,\label{cunoalfa2}
\end{eqnarray}
where $\a, 0<\a<1, C>0$ are constants depending on the \emph{a priori
data} only.

Combining the two last inequalities with \eqref{lc}, we obtain
\begin{eqnarray}\label{holder}
\|u\|_{C^{0,\a}(\G_{I})}\le
C\|u\|_{H^1(D^+_{R})}\ ,
\end{eqnarray}
where $C>0$ are constants depending on the \emph{a priori
data} only and $R=d+r_0$.
By \eqref{dipcont} we have that
 \begin{eqnarray}\label{dc2}
\|u^s\|_{H^1(D^+_{R})}\le C,
\end{eqnarray}
where $C$ is a constant depending on the \emph{a priori data} only.
Moreover, since $u=\exp{(ik\omega\cdot x)}+u^s$, by \eqref{holder}
and \eqref{dc2}, we have that
\begin{eqnarray}\label{holder1}
\|u\|_{C^{0,\a}(\G_{I})}\le C,
\end{eqnarray}
where $C$ is a constant depending on the \emph{a priori data} only.
By \eqref{holder1} and by \eqref{L}, we have that
\begin{eqnarray}
\frac{\partial u}{\partial \nu}(x)=-i\lambda(x)u(x)\ \in\
C^{0,\a}(\G_I).
\end{eqnarray}
By well-known regularity bounds for the Neumann problem (see for
instance \cite[p.667]{adn} ) it follows that, for every $R>d, \rho\in
(0,r_0),
 u\in C^{1,\a}(D^+_{R,\rho})$ and the following estimate holds
\begin{eqnarray}\label{adn}
\|u\|_{C^{1,\alpha}(D^+_{R,\rho})}\le
C_{R,\rho}\left(\|u\|_{C^{0,\alpha}(\G^{\frac{\rho}{2}}_I)}+\left\|\frac{\partial
u}{\partial
\nu}\right\|_{C^{0,\alpha}(\G^{\frac{\rho}{2}}_I)}+\|u\|_{H^1(D^+_{2R})}\right)\ ,
\end{eqnarray}
where $C_{R,\rho}>0$ is a constant depending on the \emph{a priori
data}, on $R$ and on $\rho$ only.
We shall estimate the $C^{0,\a}$ norm of $\displaystyle{\frac{\partial
u}{\partial \nu}}$ in terms of the \emph{a priori data}, indeed
\begin{eqnarray}
\left\|\frac{\partial u}{\partial
\nu}\right\|_{C^{0,\alpha}(\G^{\frac{\rho}{2}}_I))}&=& \sup_{x \in
\G^{\frac{\rho}{2}}_I}\left|\dfrac{\partial u(x)}{\partial
\nu}\right|+\left({\frac{\rho}{2}}\right)^{\a}\sup_{x,y \in
\G^{\frac{\rho}{2}}_I}\frac{\left|\dfrac{\partial u(x)}{\partial
\nu}-\dfrac{\partial
u(y)}{\partial \nu}\right|}{|x-y|^{\alpha}}=\nonumber \\
&\le&\sup_{x \in \G^{\frac{\rho}{2}}_I
}\left|\lambda(x)u(x)\right|+\left({\frac{\rho}{2}}\right)^{\a}\sup_{x,y
\in \G^{\frac{\rho}{2}}_I
}\frac{\left|\lambda(x)||u(x)-u(y)\right|}{|x-y|^{\alpha}}+\nonumber\\
&&+\left({\frac{\rho}{2}}\right)^{\a}\sup_{x,y
\in \G^{\frac{\rho}{2}}_I
}\frac{\left|u(y)||\lambda(x)-\lambda(y)\right|}{|x-y|^{\alpha}}\
.\nonumber
\end{eqnarray}
 Combining \eqref{L} and \eqref{holder1} we obtain
\begin{eqnarray}\label{derivata}
\left\|\frac{\partial u}{\partial
\nu}\right\|_{C^{0,\alpha}(\G^{\frac{\rho}{2}}_I)}&\le& \Lambda\sup_{x
\in
\G^{\frac{\rho}{2}}_I}|u(x)|+\Lambda\left({\frac{\rho}{2}}\right)^{\a}\sup_{x,y
\in
\G^{\frac{\rho}{2}}_I}\frac{|u(x)-u(y)|}{{|x-y|}^{\alpha}}+\nonumber\\
&&+\left({\frac{\rho}{2}}\right)^{\a}|\G_I|^{1-\a}\|u\|_{C^{0,\a}(\G_I)}\sup_{x,y
\in \G^{\frac{\rho}{2}}_I
}\frac{\left|\lambda(x)-\lambda(y)\right|}{|x-y|}\le\nonumber\\
 &&\le\bar{C_{\rho}}\nonumber
\end{eqnarray}
where $\bar{C_{\rho}}>0$ is a constant depending on the \emph{a priori
data} and on $\rho$
only. Moreover, since $u=\exp{(ik\omega\cdot x)}+u^s$, we have that
\eqref{dipcont} yields to
\begin{eqnarray}\label{dc1}
\|u\|_{H^1(D^+_{2R})}\le C_R,
\end{eqnarray}
where $C_R>0$ is a constant depending on the \emph{a priori data} and
on $R$ only.

Thus, inserting \eqref{holder1}, \eqref{derivata} and \eqref{dc1} in
\eqref{adn}, we obtain that
\begin{eqnarray}
\|u\|_{C^{1,\alpha}(D^+_{R,\rho})}&\le&C_{R,\rho},
\end{eqnarray}
where $C_{R,\rho}>0$ is a constant depending on the \emph{a priori
data}, on $R$ and on $\rho$ only.
\end{proof}

\begin{corollary}[\textbf{Lower bound}]\label{lb}
Let $u$ be the weak solution to \eqref{Sc}, then there exists a radius
$R_0>0$ depending on the \emph{a priori data} only, such that
\begin{eqnarray}\label{lowerbound}
|u(x)|>\frac{1}{2}\ \ \mbox{for every}\ x, \ |x|>R_0\ .
\end{eqnarray}
\end{corollary}
\begin{proof}
Let us choose  $R=4d+4r_0$. By Theorem \ref{regolarita} it follows that
there
exists a constant $C>0$ depending on the \emph{a priori data} only,
such that
\begin{eqnarray}\label{ariregolarita}
\|u\|_{C^{1,\a}\big(D^+_{2R,\frac{r_0}{2}}\big)}\le C\ .
\end{eqnarray}
In particular, by \eqref{ariregolarita}, it follows that
\begin{eqnarray}\label{ariregolarita1}
|u^s|\le C_1\ , \ \bigg|{\frac{\partial u^s}{\partial \nu}}\bigg|\le
C_1\ \ \mbox{on}\ \ \partial B_R(0),
\end{eqnarray}
where $C_1>0$ is a constant depending on the \emph{a priori data} only.

By the Green's formula for the scattered wave $u^s$ (see for instance
\cite[p.18]{ck}), we have that
\begin{eqnarray}\label{gf}
u^s(x)=\int_{\partial B_R(0)}\left(u^s(y)\frac{\partial
\phi(x,y)}{\partial
\nu(y)}-\frac{\partial u^s(y)}{\partial
\nu(y)}\phi(x,y)\right)\mbox{d}s(y), \
\ \ \ |x|>R,
\end{eqnarray}
where $$\phi(x,y)=\frac{1}{4\pi}\frac{\exp{(ik|x-y|)}}{|x-y|},\ \ x
\neq y\ ,$$
is the fundamental solution to the Helmholtz equation in $\R^3$\ .

Thus, by \eqref{gf} and by \eqref{ariregolarita1} it follows that
\begin{eqnarray}
|u^s(x)|&\le&C_1\int_{\partial B_R(0)}\bigg|\frac{\partial
\phi(x,y)}{\partial \nu(y)}\bigg|+|\phi(x,y)|\mbox{d}s(y)\le\\
&\le&C_1R^2\bigg(\frac{kR}{||x|-R|^2} + \frac{R}{||x|-R|^3} +
\frac{1}{||x|-R|}\bigg)\ .
\end{eqnarray}
Straightforward calculations show that
\begin{eqnarray}
|u^s|<\frac{1}{2}\ , \ \mbox{for every}\  x, |x|>R_0,
\end{eqnarray}
where $R_0=(k+1)8 R^3C_1 +2R$ .

The thesis follows observing that $|u|\ge 1-|u^s|$.
\end{proof}

\section{The inverse scattering problem}
\begin{lemma}[\textbf{From the far field to the near
field}]\label{farnear}
Let $u_i,u_{i,\infty},\ i=1,2$, be as in Theorem \ref{stabl}. Suppose
that, for some $\varepsilon,\ 0<\varepsilon<1$,
\eqref{farf} holds, then there exist a radius $R_1>0$ and a constant
$C>0$, depending
on the \emph{a priori data} only, such that
\begin{eqnarray}\label{nearfield}
\|u_1-u_2\|_{L^2(B_{R_1+1}(0)\setminus B_{R_1}(0))}\le C
{\varepsilon}^{\a(\varepsilon)},
\end{eqnarray}
where $\a(\varepsilon)$ is the function introduced in \eqref{lam}.
\end{lemma}
\begin{proof}
Let us choose $R=4d+4r_0$ and let us denote by $u_i^s,\ i=1,2$, the
scattered wave of the problem
\eqref{Sc} with $\lambda=\lambda_i$ respectively. By
\eqref{ariregolarita1} it follows that
\begin{eqnarray}
\|u_1^s-u_2^s\|_{L^2(\partial B_R(0))}\le C\ ,
\end{eqnarray}
where $C>0$ is a constant depending on the \emph{a priori data} only.

By the argument in \cite{is2} (see also \cite{b}), it follows that
there
exists a constant $C>0$ depending on the \emph{a priori data} only,
such that, for every $r\in (4R, 4R +1)$, the following holds
\begin{eqnarray}\label{bu}
\|u^s_1-u^s_2\|_{L^2(\partial B_r(0))}\le C
{\varepsilon}^{\a(\varepsilon)}.
\end{eqnarray}
Integrating \eqref{bu} with respect to $r$ over $(4R, 4R +1)$, we
obtain that
\begin{eqnarray}
\|u^s_1-u^s_2\|_{L^2(B_{4R +1}(0)\setminus B_{4R}(0))}\le C
{\varepsilon}^{\a(\varepsilon)}\ ,
\end{eqnarray}
 where $C>0$ is a constant depending on the \emph{a priori data} only.

Thus the thesis follows with $R_1=16d+16r_0$ and by observing that
$u^s_1-u^s_2=u_1-u_2$.

Let us stress, that H\"{o}lder stability doesn't hold, indeed, in
\cite[Section 4]{b}, it has been proved that it is not possible to
choose
$\a$ independently on $\varepsilon$. 
\end{proof}

\begin{theorem}[\textbf{Stability at the boundary}]\label{stabpc}
Let $u_i,u_{i,\infty},\ i=1,2$, be as in Theorem \ref{stabl}. We have
that there
exists $\eps_0>0$ depending on the \emph{a priori data} only, such
that,
if for some $\eps,\ 0<\eps<\eps_0$, \eqref{farf} holds, then for every
$\rho
\in (0,r_0)$ we have
\begin{eqnarray}\label{stabilitaalbordo}
\|u_1-u_2\|_{C^1(\G^{\rho}_I)}\le \eta(\eps)\ ,
\end{eqnarray}
where $\eta$ is given by \eqref{eta}, with a constant $C>0$ depending
on the \emph{a priori data} and on $\rho$ only.
\end{theorem}

\begin{proof}
By the Lipschitz regularity of the boundary $\partial D$, it follows
that the cone property holds. Namely, for every point $Q\in \partial
D$,
there exists a rigid transformation of coordinates under which we have
$Q=0$ and the finite cone
$$\mathcal{C}=\bigg\{x :|x|<r_0 ,\,\,\frac{x\cdot \xi}{|x|}>\cos\theta
\bigg \} $$
with axis in the direction $\xi$ and width $2\theta$, where
$\theta=\arctan \frac {1}{M}$, is such that $\mathcal{C}\subset D^+$.

 Let $Q$ be a point such that $Q \in \G^{r_0}_I$  and let $Q_0$ be a
point lying on the axis $\xi$ of the cone with vertex in $Q=0$ such
that
$d_0=\mbox{dist}(Q_0,0)<\frac{r_0}{2}$.

Let us define $R_2=2R_1+2$, where $R_1$ is the radius introduced in the
statement of Lemma \ref{farnear}.
Dealing as in Lieberman \cite{Li}, we consider a regularized distance
$\tilde{d}$ from the boundary of $\partial D$ such that, $\tilde{d}\in
C^2(D_{R_2}^+)\cap C^{0,1}({\overline{D_{R_2}^+}})$ and furthermore the
following properties hold\begin{itemize}
\item $\gamma_0\le \displaystyle \frac{\mbox{dist}(x,\partial
D)}{\tilde{d}(x)}\le \gamma_1$,
\item $|\nabla {\tilde{d}}(x)|\ge c_1$,\ \ \ for every $x$ such that
${\mbox {dist}}(x,\partial D)\le br_0$,
\item $\|\tilde{d}\|_{C^{0,1}}\le c_2r_0 $,
\end{itemize}
where $\gamma_0, \gamma_1, c_1, c_2, b$ are positive constants
depending on $M$ only, (see also \cite[Lemma 5.2]{abrv}).

Let us define for every $\rho>0$
\begin{eqnarray}
&& D^{\rho}=\{x\ \in D_{R_2}^+  :  \mbox{dist}(x,\partial D)>\rho \}\
,\\
&&\tilde D^{\rho}=\{x\ \in D_{R_2}^+ :  {\tilde{d}}(x)>\rho \}\ .
\end{eqnarray}
It follows that there exists $a$, $0<a\le1$, only depending on $M$
such that for every $\rho$, $0<\rho\le ar_0$, $\tilde D^{\rho}$ is
connected with boundary of class $C^1$
and
\begin{eqnarray}\label{dc}
{\tilde c_1}\rho\le \mbox{dist}(x,\partial D)\le \tilde c_2\rho \ \ \ \
\mbox{for every}\ \ x \ \in \ \partial\tilde D^{\rho},
\end{eqnarray}
where $\tilde c_1, \tilde c_2\ ,$ are positive constants depending on
$M$
only.
By\eqref{dc} we deduce that
$$D^{\tilde c_2\rho}\subset \tilde D^{\rho}\subset D^{\tilde c_1\rho}\
. $$

Let us now define ${\rho}_0=\min\{\frac{1}{16},\frac{r_0}{4}\sin
\theta\}$ and let $P$ be a point in the annulus $B_{R_1+1}(0)\setminus
B_{R_1}(0))$, such that $B_{4\rho_0}(P)\subset B_{R_1+1}(0)\setminus
B_{R_1}(0))$. Furthermore, let $\gamma$ be a path in $\tilde
D^{\frac{\rho_0}{\tilde c_1}} $ joining  $P$ to  $Q_0$ and let us
define $\{y_i\}$, $
i=0,\ldots,s$
as follows $y_0=Q_0$, $y_{i+1}=\gamma(t_{i})$, where
$t_i=\max\{\mbox{$t$ s.t. }|\gamma(t)-y_i|=2\rho_0\}$
if $|P-y_i|>2\rho_0$, otherwise let $i=s$ and stop the process.

Let us introduce the function $U\in H^1_{\mbox{loc}}(D^+)$ defined as
follows
\begin{eqnarray}\label{differenza}
U(x)=u_1(x)-u_2(x) .
\end{eqnarray}
We shall denote with $U_1$ and $U_2$ the real and the imaginary part
of $U$ respectively. Namely
$$U(x)=U_1(x)+iU_2(x) .$$
It immediately follows that $U_1,U_2,$ are both real valued solutions
to the
Helmholtz equation in $D^+$.

Thus, by the three spheres inequalities for elliptic system with
Laplacian principal part, (see \cite[Theorem 3.1]{am}), we have that
for
every $\beta_1,\beta_2,\ 1<\beta_1<\beta_2$, there exist $\bar{r}>0,
\tau,
\ 0<\tau<1$ and $C>0$ depending on the \emph{a priori data} and on
$\beta_1, \beta_2$ only, such that for every $x\in D^{\beta_2\rho}$ the
following holds
\begin{eqnarray}\label{3sfere}
 \int_{B_{\beta_1\rho}(x)} |U|^2 \le C \left(\int_{B_{\rho}(x)}
|U|^2\right)^{\tau}\cdot \left(\int_{B_{\beta_2\rho}(x)}
|U|^2\right)^{1-\tau}
\end{eqnarray}
for every $\rho\in (0,\bar {r})$.
By a possible replacement of $\rho_0$ with $\bar{r}$
if $\rho_0>\bar{r}$ and choosing in \eqref{3sfere} $\beta_1=3,\
\beta_2=4, \rho=\rho_0, \
x=y_0$, we infer that
\begin{eqnarray}\label{ddddd}
 \int_{B_{3\rho_0}(y_0)} |U|^2 \le C\left(\int_{B_{\rho_0}(y_0)}
|U|^2\right)^{\tau}\cdot \left(\int_{B_{4\rho_0}(y_0)}
|U|^2\right)^{1-\tau}.
\end{eqnarray}

As a consequence of Lemma \ref{benposto}, we have that
\begin{eqnarray}\label{diffbenposto}
 \|U\|_{H^1(D_{R_2}^+)}\le C,
\end{eqnarray}
where $C>0$ is a constant depending on the \emph{a priori data} only.

Let us observe that $B_{4\rho_0}(y_0)\subset D_{R_2}^+$ and
$B_{\rho_0}(y_0)\subset
B_{3\rho_0}(y_1)$. Thus by \eqref{ddddd} and \eqref{diffbenposto}
we deduce that
$$ \int_{B_{\rho_0}(y_0)}|U| ^2 \le C\left(\int_{B_{3\rho_0}(y_1)}
|U|^2\right)^{\tau}\cdot C^{1-\tau}\ .$$
An iterated application of the three spheres inequality leads to
$$  \int_{B_{\rho_0}(y_0)} |U|^2 \le   \left(\int_{B_{\rho_0}(y_s)}
|U|^2\right)^{{\tau}^{s}}\cdot C^{1-{\tau}^s}\ .$$
Finally, since $B_{\rho}(y_s)\subset B_{R_1+1}(0)\setminus
B_{R_1}(0))$,  by \eqref{nearfield} we obtain
that
$$ \int_{B_{\rho_0}(y_0)} |U|^2 \le C\big\{
\eps^{\a{(\eps)}}\big\}^{{\tau}^s}\ .$$
We shall construct a chain of balls $B_{\rho_k}(Q_k)$ centered on the
axis of the cone, pairwise tangent to each other and all contained in
the cone
$${\mathcal{C}}^{\prime}=\bigg\{x :|x|<r_0 ,\,\,\frac{x\cdot
\xi}{|x|}>\cos{\theta}^{\prime} \bigg \}\ ,$$
where $\theta^{\prime}=\arcsin\big(\frac{\rho_0}{d_0}\big) .$
Let $B_{\rho_0}(Q_0)$ be the first of them, the following are defined
by induction in such a way
\begin{equation*}
\begin{array}{l}
Q_{k+1}=Q_{k}-(1+\mu)\rho_k\xi\ ,
\\
\rho_{k+1}=\mu\rho_{k}\ ,
\\
d_{k+1}=\mu d_{k}\ ,
\end{array}
\end{equation*}
with
\begin{equation*}
\begin{array}{l}
 \mu=\dfrac{1-\sin\theta^{\prime}}{1+\sin\theta^{\prime}}\ .
\end{array}
\end{equation*}
Hence, with this choice, we have $\rho_{k}=\mu^{k}\rho_0$ and
$B_{\rho_{k+1}}(Q_{k+1})\subset B_{3\rho_k}(Q_k)$.\\
Considering the following estimate obtained by a repeated application
of the three spheres inequality, we have that
\begin{eqnarray}\label{ss}
\| U \|_{L^2(B_{\rho_k}(Q_k))}&\le& \| U
\|_{L^2(B_{3\rho_{k-1}}(Q_{k-1}))
}\le \nonumber\\
&\le&  \| U \|^{\tau}_{L^2(B_{\rho_{k-1}}(Q_{k-1}))}  \| U
\|^{1-\tau}_{L^2(B_{4\rho_{l-1}}(Q_{k-1}))}\nonumber\\
&\le& C \| U  \|^{{\tau}^k}_{L^2(B_{{\rho}_0}(Q_0))}\le C\Big\{ \big[
\eps^{\a(\eps)}\big]^{{\tau}^s}\Big\}^{{\tau}^k}\ .
\end{eqnarray}
For every $r$, $0<r<d_0$, let $k(r)$ be the smallest positive integer
such that  $d_k\le r$ then, since $d_k={\mu}^k d_0$, it follows
\begin{eqnarray}\label{r}
\dfrac{|\log(\frac{r}{d_0})|}{\log{\mu}}\le k(r)\le
\dfrac{|\log(\frac{r}{d_0})|}{\log{\mu}} +1\ ,
\end{eqnarray}
and by \eqref{ss} we deduce
\begin{eqnarray}\label{k}
   \| U \|_{L^2(B_{\rho_k(r)}(Q_k(r)))}\le
  C  \Big\{ \big[
\eps^{\a(\eps)}\big]^{{\tau}^s}\Big\}^{{\tau}^{k(r)}}\ .
\end{eqnarray}

Let $\bar{x} \in \G^{\frac{\rho}{2}}_I$ with $\rho\in (0,r_0)$ and let
$x \in B_{\frac{\rho_{k(r)-1}}{2}}(Q_{k(r)-1}) $. By Theorem
\ref{regolarita}, in particular, it follows that $U \in
C^{1,\a}(D^+_{R_2,
\frac{\rho}{4}})$ with
\begin{eqnarray}\label{ciunoalfaU}
\| U\|_{C^{1,\a}(D^+_{R_2, \frac{\rho}{4}})}\le C_{\rho},
\end{eqnarray}
where $C_{\rho}>0$ is a constant depending on the \emph{a priori data}
and on $\rho$ only.
 Then \eqref{ciunoalfaU} yields to
$$\left|U(\bar{x})\right| \le
 \left|U({x})\right| +C_{\rho}|x-\bar{x}|^{\a} \le
\left|U({x})\right| +C_{\rho}\bigg( {\frac{2}{\mu} r}\bigg)^{\a}\ . $$

Integrating this inequality over $
B_{\frac{\rho_{k(r)-1}}{2}}(Q_{k(r)-1}) $, we have that

\begin{eqnarray}\label{integrato}
\left| U(\bar{x})\right|^2&\le&\!\!\!
\frac{2}{{\omega_3}{(\frac{\rho_{k-1}}{2})}^3}\int_{B_{\frac{\rho_{k(r)-1}}{2}}\big(Q_{k(r)-1}\big)}
| U(x)|^2 \mbox{d}x + 2 C_{\rho}^2{\bigg(\frac{4
r^2}{{\mu}^2}\bigg)}^{\a}\!\!\!.
\end{eqnarray}

Being $k$ the smallest integer such that $d_k\le r$, then $d_{k-1}>r$
and thus \eqref{integrato} yields to
$$\left| U(\bar{x})\right|^2\le\frac{C}{\big({r \sin
{\theta}^{\prime}\big)}
^{3}}  \int_{ B_{\rho_{k(r)-1}}(Q_{k(r)-1})} |U(x)|^2  \mbox{d}x +
C_{\rho} r^{2\a}\ .$$

By  $\eqref {k}$ we deduce that
\begin{eqnarray}\label{funz}
\left|U(\bar{x})\right|^2\le
\frac{C}{{r}^{3}} \Big\{ \big[\eps^{\a{(\eps)}}\big]^{{\tau}^s}
\Big\}^
{{\tau}^{k(r)-1}} +C_{\rho}r^{2\a}\ .
\end{eqnarray}
The estimate \eqref{ciunoalfaU} also provides us that
$$\left|\frac{\partial U(\bar{x})}{\partial \nu}\right| \le
\left|\dfrac{\partial U({x})}{\partial \nu}\right| +C_{\rho}\bigg(
{\frac{2}{\mu} r}\bigg)^{\a}\ .$$
Integrating over $ B_{\frac{\rho_{k(r)-1}}{2}}(Q_{k(r)-1}) $ we deduce
that
\begin{eqnarray}
\left|\frac{\partial U(\bar{x})}{\partial \nu}\right|^2&\le&
\frac{2}{{\omega_3}{(\frac{\rho_{k-1}}{2})}^3}\int_{B_{\frac{\rho_{k(r)-1}}{2}}\big(Q_{k(r)-1}\big)}
\bigg| \frac{\partial U(x)}{\partial \nu} \bigg|^2 \mbox{d}x + 2
C_{\rho}^2{\bigg(\frac{4 r^2}{{\mu}^2}\bigg)}^{\a}\le \nonumber\\
&\le&
\frac{2}{{\omega_3}{(\frac{\rho_{k-1}}{2})}^3}\int_{B_{\frac{\rho_{k(r)-1}}{2}}\big(Q_{k(r)-1}\big)}
|\nabla U(x)|^2 \mbox{d}x + 2 C_{\rho}^2{\bigg(\frac{4
r^2}{{\mu}^2}\bigg)}^{\a}\ .\nonumber
\end{eqnarray}
Applying the Caccioppoli inequality, we have
$$\left|\frac{\partial U(\bar{x})}{\partial \nu}\right|^2\le
\frac{C}{{\big(\rho_{k-1}\big)}
^{5}}\int_{ B_{\rho_{k(r)-1}}(Q_{k(r)-1})} U(x)^2 \mbox{d}x +C_{\rho}
r^{2\a} \ .$$

Dealing with the same arguments that lead to \eqref{funz}, we obtain
that
\begin{eqnarray}\label{dernor}
\left|\frac{\partial U(\bar{x})}{\partial \nu}\right|^2\le
\frac{C}{{r}^{5}} \Big\{ \big[\eps^{\a{(\eps)}}
\big]^{{\tau}^s}\Big\}^{{\tau}^{k(r)-1}} +C_{\rho}r^{2\a}\ .
\end{eqnarray}

The choice in $\eqref{r}$ guarantees that
$$\tau^{k(r)-1}\ge \bigg(\frac{r}{d_0}\bigg)^{\nu}, $$
where $\nu= -\log\big(\frac{1}{\mu}\big)\log\tau$.
Thus, by \eqref{funz} and by \eqref{dernor}, it follows that
\begin{eqnarray}
&& \left| U(\bar{x})\right|\le C_{\rho}
\bigg\{r^{-{\frac{3}{2}}}\Big[\big(\eps^{\a{(\eps)}}
\big)^{{\tau}^s}\Big]^{\frac{r^{\nu}}{2}}+r^{\a}
\bigg\}\ ,\\
&& \left|\frac{\partial U(\bar{x})}{\partial \nu}\right|\le C_{\rho}
\bigg\{r^{-{\frac{5}{2}}}\Big[\big(\eps^{\a{(\eps)}}
\big)^{{\tau}^s}\Big]^{\frac{r^{\nu}}{2}}+r^{\a} \bigg\}\ .
\end{eqnarray}
Minimizing the right hand sides of the above inequalities with respect
to $r$, with $r \in (0,\frac{r_0}{4})$, we deduce
\begin{eqnarray}
&&\left| U(\bar{x})\right|\le C_{\rho}\big(\log{(\eps^{-\a(\eps)})}
\big)
^{-\frac{2\a}{\nu+2}}\ ,\label{duestime1}\\
&&\left|\frac{\partial U(\bar{x})}{\partial \nu}\right|\le
C_{\rho}\big(\log{(\eps^{-\a(\eps)})}
\big)^{-\frac{2\a}{\nu+2}},\label{duestime2}
\end{eqnarray}
where $C_{\rho}>0$ is a constant depending on the \emph{a priori data}
and on $\rho$ only.
Thus, since $\bar{x}$ is an arbitrary point in $\G_I^{\frac{\rho}{2}}$,
by \eqref{duestime1} and \eqref{duestime2} we have that
\begin{eqnarray}\label{duestimeinf}
&&\| U(\bar{x})\|_{L^{\infty}(\G_I^{\frac{\rho}{2}})}\le
C_{\rho}\big(\log{(\eps^{-\a(\eps)})} \big)
^{-\frac{2\a}{\nu+2}}\ ,\\
&&\bigg\|\frac{\partial U(\bar{x})}{\partial
\nu}\bigg\|_{L^{\infty}(\G_I^{\frac{\rho}{2}})}\le
C_{\rho}\big(\log{(\eps^{-\a(\eps)})}
\big)^{-\frac{2\a}{\nu+2}}.
\end{eqnarray}

By an interpolation inequality we have
$$\|\nabla_t(U)\|_{ L^{\infty}(\G_{1,\rho})}\le c_{\rho} \|U\|_{
L^{\infty}(\G_{1,\frac{\rho}{2}})} ^{\beta}{\|U
\|_{C^{1,\a}(\G_{1,\rho})}}^{1-\beta}\ ,$$
where $\beta=\frac{\a}{\a+1}$ and $c_{\rho}>0$ depends on the \emph{a
priori data} and on $\rho$ only.
Thus, by \eqref{ciunoalfaU}, we obtain
$$\|\nabla_t(U)\|_{ L^{\infty}(\G_{1,\rho})}\le c_{\rho}\|U\|_{
L^{\infty}(\G_{1,\frac{\rho}{2}})} ^{\beta}{C_{\rho}}^{1-\beta}\ .$$

It follows that for every $\eps<\eps_0$, with $\eps_0$ depending only
on the \emph{a priori data},
\begin{eqnarray}
\|\nabla(U)\|_{ L^{\infty}(\G_{1,\rho})}&\le&\left\|
\frac{\partial{U}}{\partial{\nu}}\right \|
_{L^{\infty}(\G_{1,\rho})} +
\|\nabla_t(U)\|_{ L^{\infty}( \G_{1,\rho})}\le\nonumber
\\
&\leq& {C}_{\rho}\big(\log{(\eps^{-\a(\eps)})} \big)
^{-\frac{2\a \beta}{\nu+2}},
\end{eqnarray}
where $C_{\rho}>0$ depends on the \emph{a priori data} and on $\rho$
only.
Hence, by a possible replacing of $\eps_0$ with a smaller one depending
on the \emph{a priori data} only, we have that
\begin{eqnarray}
\|u_1-u_2 \|_{ C^1(\G_{1,\rho})}\leq
{C}_{\rho}\big(\log{(\eps^{-\a(\eps)})}\big)^{-\frac{2\a
\beta}{\nu+2}}\ \ \mbox{for every}\ \eps,
0<\eps<\eps_0.
\end{eqnarray}
Thus the thesis follows replacing in \eqref{eta} $C$ with $C_{\rho}$
and
$\theta$ with $\frac{2\a \beta}{\nu+2}$.
\end{proof}
\begin{proposition}\label{psi2}
There exists a radius $r_1>0$ depending on the \emph{a priori data}
only such that, for every $x_0\in\G^{r_0}_I$, the problem
\begin{equation}\label{psiz}
\left\{
\begin{array}
{lcl}
\Delta \psi +k^2\psi=0, &\   \mbox{in $\G_{I,r_1}(x_0)$},
\\
\dfrac{\partial \psi}{\partial \nu} + i\lambda(x)\psi=0, &\ \mbox{on
$\Delta_{I,r_1}(x_0)$},
\end{array}
\right.
\end{equation}
admits a solution $\psi\in H^1(\G_{I,r_1}(x_0))$ satisfying
\begin{eqnarray}\label{uno}
|\psi(x)|\ge 1 \ \mbox{for every}\  x\in\G_{I,r_1}(x_0).
\end{eqnarray}
Moreover, there exists a constant $\bar{\psi}>0$ depending on the
\emph{a
priori data} only, such that for every $x_0\in\G^{r_0}_I$
\begin{eqnarray}\label{psize}
\|\psi\|_{C^{1}(\G_{I,r_1}(x_0))}\le {\bar{\psi}}.
\end{eqnarray}
\end{proposition}
\begin{proof}
Let us consider a point $x_0\in \G^{r_0}_I$. After a translation we may
assume that $x_0=0$ and, fixing local coordinates, we can represent the
boundary as a graph of a $C^{1,1}$ function. Namely, we have that
\begin{eqnarray}\label{loccoor}
D^+ \cap B_{r_0}(0)=\{(x',x_3) \in  B_{r_0}(0): x_3<\varphi_I(x')\}\ ,
\end{eqnarray}
where $\varphi_I$ is the $C^{1,1}$ function satisfying
\eqref{sI1},\eqref{sI2},\eqref{sI3}.

Let $\Phi\in C^{1,1}(B_{\frac{r_0}{4M}},\R^3)$ be the map defined as
follows

\begin{eqnarray}
\Phi(y',y_3)=(y',y_3+\varphi_I(y'))\ .
\end{eqnarray}
We have that there exist $\theta_1,\theta_2,\theta_1>1>\theta_2>0,$
constants depending on $M$ and $r_0$ only, such that, for every $r\in
(0,\frac{r_0}{4M})$, it follows that
\begin{eqnarray}\label{inclusioni}
\G_{I, \theta_2 r}(0)\subset \Phi(B^-_r(0))\subset \G_{I, \theta_1
r}(0)\ ,
\end{eqnarray}
where $B^-_r(0)=\{y\in\R^3: |y|<r,\ y_3<0 \}$ and furthermore we have
\begin{eqnarray}
|\mbox{det}\ D\Phi\ |=1\ .
\end{eqnarray}
The inverse map $\Phi^{-1}\in C^{1,1}(\G_{I,r_0}(0),\R^3)$ and is
defined
by
\begin{eqnarray}
\Phi^{-1}(x',x_3)=(x',x_3-\varphi_I(x'))\ .
\end{eqnarray}

Denoting by
\begin{eqnarray}\label{nuovacondutt}
\sigma(y)=(\sigma_{i,j}(y))_{i,j=1}^3=(D\Phi^{-1})(\Phi(y))\cdot(D\Phi^{-1})^T(\Phi(y))\
,
\end{eqnarray}
\begin{eqnarray}\label{nuovaimpedenza}
\lambda'(y)=\lambda(\Phi(y))\ ,
\end{eqnarray}
\begin{eqnarray}\label{lambdazero}
{\lambda_0}'=\lambda'(0)\ ,
\end{eqnarray}

it follows that
\begin{eqnarray}\label{inzero}
\sigma(0)=\mbox{I},
\end{eqnarray}
\begin{eqnarray}\label{reg}
\ \ \ \|\sigma_{i,j}\|_{C^{0,1}(\G_{I,r_0 })}\le \Sigma,\ \ \
\mbox{for}\ \ i,j=1,2,3,
\end{eqnarray}

\begin{eqnarray}\label{unifell}
\frac{1}{2}|\xi|^2\le\sigma(y)\xi\cdot{\xi}\le C_1|\xi|^2, \ \mbox{for
every}\ y\in B^-_{(\frac{r_0}{4M})}(0)\  \mbox{and every}\
\xi\in\mathbb{R}^3,
\end{eqnarray}
\begin{eqnarray}\label{regl}
\|\lambda'\|_{C^{0,1}(B'_{\frac{r_0}{4M}}(0))}\le \Lambda'\ ,
\end{eqnarray}
where $\Sigma>0, C_1>0, \Lambda'>0$ are constants depending on
$M,r_0,\Lambda$ only.

\begin{claim}\label{cl}
There exists a radius ${r_2},\ 0<{r_2}<\frac{r_0}{4M}$ and a solution
$\psi'\in H^1(B^-_{r_2}(0))$ to the problem
\begin{equation}\label{psip}
\left\{
\begin{array}
{lcl}
\mbox{div}(\sigma\nabla {\psi}') +k^2{\psi}'=0\ ,&\ \ \  \mbox{in}\ \ \
B^-_{r_2}(0)\ ,
\\
\sigma\nabla {\psi}'\cdot\nu' + i\lambda' {\psi}'=0\ , &\ \mbox{on}\ \
\ B'_{r_2}(0),
\end{array}
\right.
\end{equation}
where $\nu'=(0,0,1)$ such that
$$|\psi'|\ge 1\  \mbox{in} \ B^-_{r_2}(0).$$
\end{claim}




\begin{proof}\!\!\!\!\textbf{of Claim \ref{cl}.}

We look for a radius $r_2>0$ and for a solution of the form
${\psi}'=\psi_0 -s$ such that, $\psi_0
\in H^1(B^-_{r_2}(0))$ is a weak solution to the problem
\begin{equation}\label{psizero}
\left\{
\begin{array}
{lll}
\Delta \psi_0 +k^2\psi_0=0\ ,& \mbox{in}\ \ \ B^-_{r_2}(0),
\\
\dfrac{\partial \psi_0}{\partial \nu} + i{\lambda_0}'\psi_0=0\ , \ &
\mbox{on}\ \ \ B'_{r_2}(0),
\end{array}
\right.
\end{equation}
satisfying $|\psi_0|\ge 2$ in $B^-_{r_2}(0)$ .

Whereas $s\in H^1(B^-_{r_2}(0))$ is a weak solution to the problem
\begin{equation}\label{s}
\left\{
\begin{array}
{lll}
\mbox{div}(\sigma\nabla s) +k^2s=\mbox{div}((\sigma-I)\nabla \psi_0)\ ,
&\mbox{in}\ \ \ B^-_{r_2}(0),
\\
\sigma\nabla s\cdot\nu + i\lambda' s=(\sigma-I)\nabla \psi_0\cdot\nu +
i(\lambda'-{\lambda_0}')\psi_0\ ,  &\mbox{on}\ \ \ B'_{r_2}(0) ,
\\
s=0\ , \ \ \ \ &\mbox{on}\ \ \ |y|={r_2},
\end{array}
\right.
\end{equation}
such that $s(y)=O(|y|^{2})$ near the origin.

We can construct $\psi_0$ explicitly as follows
\begin{eqnarray*}
&&\psi_0(y_1,y_2,y_3)=8\cosh\big(|{\lambda_0}'^2-k^2|^{\frac{1}{2}}y_1\big)\big[\sin\big({\lambda_0}'
y_3 \big)+i\cos\big({\lambda_0}'
y_3\big)\big],\ \mbox{if}\ k^2<{\lambda_0}'^2,\\
&&\psi_0(y_1,y_2,y_3)=8\cos\big(|k^2-{\lambda_0}'^2|^{\frac{1}{2}}y_1\big)\big[\sin\big({\lambda_0}'
y_3\big)+i\cos\big({\lambda_0}'
y_3\big)\big],\ \mbox{if}\ k^2>{\lambda_0}'^2,\\
&&\psi_0(y_1,y_2,y_3)=8\sin\big({\lambda_0}'
y_3\big)+i8\cos\big({\lambda_0}' y_3\big),\
\mbox{if}\ k^2={\lambda_0}'^2.
\end{eqnarray*}
Denoting by
\begin{eqnarray}
\tilde{r}=\frac{\pi}{4}\min\left\{\frac{1}{\sqrt{|k^2-{\lambda_0}'^2|}},
\frac{1}{{\lambda_0}'}\right\}\ ,
\end{eqnarray}
it follows, by straightforward calculations, that $\psi_0\in
H^1(B^-_{\tilde{r}}(0))$ is a weak solution of \eqref{psizero}
with $r_2=\tilde{r}$
and $|\psi_0|\ge 2$ in $B^-_{\tilde{r}}(0)$.

Let us now look for a solution $s$ to the problem \eqref{s}.

Fixed $r\in(0,\frac{r_0}{8M})$, let us define the space
\begin{eqnarray}\label{spazio}
H^1_{0^-}(B^-_r(0))=\{\eta\in H^1(B^-_r(0))\  \mbox{such that}\
\eta(y)=0\ \mbox{on}\ \ |y|=r\},
\end{eqnarray}
endowed with the usual $\|\cdot\|_{H^1_{0}(B^-_r(0))}$ norm.
Thus the weak formulation of the problem \eqref{s} reads in this way:
find $s\in H^1_{0^-}(B^-_r(0))$ such that, for every $\eta \in
H^1_{0^-}(B^-_r(0))$, the following holds

\begin{eqnarray}\label{wfr}
\int_{B^-_r(0)}\sigma\nabla s\cdot\nabla\bar{\eta}
-\int_{B^-_r(0)}k^2s\bar{\eta}-\int_{B'_r(0)}i\lambda'
s\bar{\eta}&&=\int_{B^-_r(0)}
(\sigma -I)\nabla \psi_0\cdot\nabla \bar{\eta} +\nonumber\\
&&+i\int_{B'_r(0)}(\lambda'-{\lambda_0}') \psi_0\bar{\eta}.
\end{eqnarray}
Let us introduce the following bilinear form
\begin{eqnarray}
A&:&H^1_{0^-}(B^-_r(0)) \times H^1_{0^-}(B^-_r(0)) \rightarrow
\mathbb{C}
\end{eqnarray}
such that
\begin{eqnarray}
A(\eta_1,\eta_2)=\int_{B^-_r(0)}\sigma\nabla
\eta_1\cdot\nabla\bar{\eta_2}
-\int_{B^-_r(0)}k^2\eta_1\bar{\eta_2}-\int_{B'_r(0)}i\lambda'
\eta_1\bar{\eta_2}
\end{eqnarray}
and  the following functional
\begin{eqnarray}
F:H^1_{0^-}(B^-_r(0))\rightarrow \mathbb{C}
\end{eqnarray}
such that
\begin{eqnarray}
F(\eta)=\int_{B^-_r(0)}
(\sigma -I)\nabla \psi_0\cdot\nabla
\bar{\eta}+i\int_{B'_r(0)}(\lambda'-{\lambda_0}') \psi_0\bar{\eta}\ .
\end{eqnarray}
It immediately follows that $A$ and $F$ are continuous on
$H^1_{0^-}(B^-_r(0))$ as bilinear form and as a functional
respectively.

Moreover, dealing as in \cite[Lemma 8.4]{gt}, we have that, by the
H\"{older} inequality, it follows that for every $\eta
\in H^1_{0^-}(B^-_r(0))$
\begin{eqnarray}\label{poincare2}
 \int_{B^-_r(0)}|\eta|^2\le \tilde{c_1}r^2
\Big(\int_{B^-_r(0)}|\eta|^6\Big)^{\frac{1}{3}},
\end{eqnarray}
where $\tilde{c_1}>0$ is a constant depending on the \emph{a priori
data} only.
Hence by the Sobolev Imbedding Theorem, (see \cite[Chap.4]{Adams}), and
by \eqref{poincare2}, we have that
\begin{eqnarray}\label{poincare}
 \int_{B^-_r(0)}|\eta|^2\le c_1r^2 \int_{B^-_r(0)}|\nabla\eta|^2,
\end{eqnarray}
where $c_1>0$ is a constant depending on the \emph{a priori data} only.

Analogously, by the H\"{older} inequality on the boundary, it follows
that
\begin{eqnarray}\label{traccia1}
 \int_{B'_r(0)}|\eta|^2\le
\tilde{c_2}r\Big(\int_{B'_r(0)}|\eta|^4\Big)^{\frac{1}{2}},
\end{eqnarray}
where $\tilde{c_2}>0$ is a constant depending on the \emph{a priori
data} only.
By a trace inequality (see for instance \cite{Adams}, Chap.5), it
follows that
\begin{eqnarray}\label{traccia}
 \int_{B'_r(0)}|\eta|^2\le c_2r\int_{B^-_r(0)}|\nabla\eta|^2,
\end{eqnarray}
where $c_2>0$ is a constant depending on the \emph{a priori data} only.

 Thus, by \eqref{unifell},\eqref{poincare} and \eqref{traccia}, we
deduce that
\begin{eqnarray*}
|A(\eta,\eta)|\ge( \frac{1}{2}-
c_1r^2k^2-c_2r\Lambda')\int_{B^-_r(0)}|\nabla \eta|^2 .
\end{eqnarray*}
Denoting by
\begin{eqnarray}
r_3=\min\left\{1,\frac{1}{8}(c_1k^2+c_2\Lambda),\frac{r_0}{8M}\right\},
\end{eqnarray}
we have that for every $r\in (0,r_3)$
\begin{eqnarray}\label{coercitivita}
|A(\eta,\eta)|\ge\frac{1}{4}\int_{B^-_r(0)}|\nabla \eta|^2.
\end{eqnarray}
Thus it follows that, for every $r\in(0,r_3)$, the bilinear form $A$ is
coercive on
$H^1_{0^-}(B^-_{r}(0))$. Hence by the Lax-Milgram theorem we can infer
that, for every $r\in(0,r_3)$, there exists a unique solution $s\in
H^1_{0^-}(B^-_{r}(0))$ to
the
problem \eqref{s}.

Fixing $r\in(0,r_3)$ and choosing $\eta=s$ as test function in the weak
formulation \eqref{wfr},
we obtain
\begin{eqnarray}\label{rconiug}
\int_{B^-_r(0)}\sigma\nabla s\cdot\nabla\bar{s}
-\int_{B^-_r(0)}k^2|s|^2-\int_{B'_r(0)}i\lambda'|s|^2&=&\int_{B^-_r(0)}
(\sigma -I)\nabla \psi_0\cdot\nabla \bar{s} +\nonumber\\
+&&\!\!\!\!\!i\int_{B'_r(0)}(\lambda'-{\lambda_0}') \psi_0\bar{s}.
\end{eqnarray}
By \eqref{coercitivita}, we have that
\begin{eqnarray}\label{x}
\frac{1}{4}\int_{B^-_r(0)}|\nabla s|^2\le \Big|\int_{B^-_r(0)}
(\sigma -I)\nabla \psi_0\cdot\nabla
\bar{s}\Big|+\Big|\int_{B'_r(0)}(\lambda'-{\lambda_0}')
\psi_0\bar{s}\Big|.
\end{eqnarray}

By the Schwartz inequality, by \eqref{inzero} and by \eqref{reg} we
have that
\begin{eqnarray}\label{y}
\Big|\int_{B^-_r(0)}
(\sigma -I)\nabla \psi_0\cdot\nabla \bar{s}\Big|\le 16\Sigma
r^2\int_{B^-_r(0)}|\nabla \psi_0|^2 +
\frac{1}{16}\int_{B^-_r(0)}|\nabla s|^2\
.\end{eqnarray}
Analogously, we have that, by the Schwartz inequality, by
\eqref{lambdazero}
and by \eqref{regl} it follows that
\begin{eqnarray}\label{yy}
\Big|\int_{B'_r(0)}(\lambda'-{\lambda_0}') \psi_0\bar{s}\Big|\le
16c_2\Lambda'r^2\int_{B'_r(0)}|\psi_0|^2+\frac{1}{16c_2}\int_{B'_r(0)}|s|^2\
.
\end{eqnarray}
Moreover, by the inequality \eqref{traccia} and by \eqref{yy} we deduce
\begin{eqnarray}\label{yyy}
\Big|\int_{B'_r(0)}(\lambda'-{\lambda_0}') \psi_0\bar{s}\Big|\le c_2^2
r^4
16\Lambda'\int_{B^-_r(0)}|\nabla
\psi_0|^2+\frac{1}{16}r\int_{B^-_r(0)}|\nabla s|^2\ .
\end{eqnarray}
Hence inserting \eqref{y} and \eqref{yyy} in \eqref{x} we obtain that
\begin{eqnarray}
\frac{1}{8}\int_{B^-_r(0)}|\nabla s|^2\le (16\Sigma +c_2^2
16\Lambda')r^2\int_{B^-_r(0)}|\nabla \psi_0|^2\ .
\end{eqnarray}

Denoting by
 $$Q=\sup_{B^-_{\frac{r_0}{8M}}(0)}|\nabla \psi_0|^2 ,$$
we have that
\begin{eqnarray}
\frac{1}{8}\int_{B^-_r(0)}|\nabla s|^2\le \frac{4}{3}\pi(16\Sigma
+c_1^2 16\Lambda')r^5Q\ .
\end{eqnarray}
By standard estimates for solutions of elliptic equations (see for
instance \cite{gt}, Chap.8) and observing that $Q>0$ depends on the
\emph{a priori data} only, we can infer that for every
$r\in(0,\frac{r_3}{2})$
\begin{eqnarray*}
\|s\|_{L^{\infty}(B^-_{r}(0))}\le c_4 r^2,
\end{eqnarray*}
where $c_4>0$ is a constant depending on the \emph{a priori data} only.

Hence the Claim follows choosing ${r_2}=\min\{\tilde{r},\frac{r_3}{2},
\frac{1}{\sqrt{c_4}}\}$ and observing that
$$|{\psi}'|\ge|\psi_0|-|s|\ge 1\ \ \ \mbox{in}\ B^-_{{r_2}}(0)\ .$$
\end{proof}
Let us notice that choosing $r_1=\theta_2{r_2}$ and
$\psi(x',x_3)=\psi'(\Phi^{-1}(x',x_3))$, we have that $\psi\in
H^1(\G_{I,r_1}(0))$ is a weak
solution to the problem \eqref{psiz} and is such that $|\psi|\ge 1$ in
$\G_{I,r_1}(0)$.

Finally, we conclude the proof of Proposition \ref{psi2} observing that
\eqref{psize} follows
dealing with
the same argument used in the proof of Theorem \ref{regolarita}.
\end{proof}

\begin{lemma}[Volume doubling inequality]\label{dib}
Let $u$ be the solution to the problem \eqref{Sc}, then there exists a
radius
$\bar{r}>0$ such that for every $x_0\in
\G^{r_0}_{I}$ the following holds
\begin{eqnarray}\label{doublingat}
\int_{\Gamma_{I, \beta r}}|u|^2 \le C{\b}^K\int_{\Gamma_{I,r}}|u|^2
 \end{eqnarray}
for every $r,\beta$ such that $\b>1$ and $0<\b r<\bar{r}$, where $C>0,
K>0$ are constants depending on the \emph{a priori data} only.
\end{lemma}
\begin{proof}
Let $x_0\in\G^{r_0}_{I}$ and let $r_1$ and $\psi$ be, respectively, the
radius and the
function, introduced in Proposition \ref{psi2}. Denoting by
\begin{eqnarray}\label{defz}
z=\frac{u}{\psi},
\end{eqnarray}
 it follows that $z\in H^1(\G_{I,r_1}(x_0))$ is a weak solution to the
problem
\begin{equation}\label{z2}
\left\{
\begin{array}
{lcl}
\Delta z + 2\displaystyle{\frac{\nabla \psi}{\psi}}\cdot\nabla z=0,& \
\ \mbox{in
$\G_{I,r_1}(x_0)$},
\\
\dfrac{\partial z}{\partial \nu}=0,&\ \  \mbox{on
$\Delta_{I,r_1}(x_0)$}.
\end{array}
\right.
\end{equation}
Dealing as in Proposition \ref{psi2}, we may assume that, up to a rigid
transformation of coordinates, $x_0=0$ and, by local coordinates, we
can
locally represent the boundary as a graph of a $C^{1,1}$ function as in
\eqref{loccoor}.

Following \cite[Theorem 0.8]{ae}, (see also \cite[Proposition
3.5]{abrv}), we have that there exists a map $\Psi\in
C^{1,1}(B_{\rho_2}(0),\R^3)$ such that
\begin{eqnarray}
\Psi(B_{\rho_2}(0))\subset B_{\rho_1}(0),
\end{eqnarray}
\begin{eqnarray}
\Psi(y',0)=(y',\varphi_I(y')),\ \ \ \ \ \ \mbox{for every}\ y'\in
B'_{\rho_2}(0),
\end{eqnarray}
\begin{eqnarray}\label{inclusioni2}
\ \ \ \G_{I,\frac{\rho}{2}}\subset
\Psi(B^-_{\rho}(0))\subset\G_{I,c_1\rho},\ \ \mbox{for every}\ \rho\in
(0,\rho_2),
\end{eqnarray}
\begin{eqnarray}
\frac{1}{8}\le |\mbox{det} D\Psi|\le c_2,
\end{eqnarray}
where $\rho_1, 0<\rho_1<r_0,\rho_2>0,c_1>0,c_2>0$ are constants
depending on $r_0, M, \Lambda$ only.
Denoting by
\begin{eqnarray}
A(y)=|\mbox{det}
D\Psi(y)|(D\Psi^{-1})(\Psi(y))(D\Psi^{-1})^T(\Psi(y)),
\end{eqnarray}
\begin{eqnarray}
B(y)=2|\mbox{det}
D\Psi(y)|(D\Psi^{-1})(\Psi(y))\frac{\nabla
\psi(\Psi(y))}{\psi(\Psi(y))},
\end{eqnarray}
\begin{eqnarray}
v(y)=z(\Psi(y)),
\end{eqnarray}
it follows that
\begin{eqnarray}
A(0)=I\ ,
\end{eqnarray}
\begin{eqnarray}
A(y',0)(y',0)\cdot e_3=0,\ \mbox{for every}\ \ y', |y'|\le \rho_2,
\end{eqnarray}
\begin{eqnarray}\label{adoesc}
c_3|\xi|^2\le A(y)\xi\cdot\xi\le c_4|\xi|^2, \mbox{for every}\ \ y\in
B^-_{\rho_2}(0)\ \mbox{and for every}\ \xi\in\mathbb{R}^3,
\end{eqnarray}
\begin{eqnarray}
|A(y_1)-A(y_2)|\le c_5|y_1-y_2|,\ \ \mbox{for every}\ \ y_1,y_2\in
B^-_{\rho_2}(0),
\end{eqnarray}
\begin{eqnarray}
|B(y)|\le c_6,\ \mbox{for every}\ \ y\in
B^-_{\rho_2}(0),
\end{eqnarray}
where $c_4>0,c_5>0,c_6>0$ are constants depending on $r_0, M, \Lambda$
only.

Let us observe that $v\in H^1(B^-_{\rho_2}(0))$ is a weak solution to
the problem
\begin{equation}\label{v}
\left\{
\begin{array}
{lcl}
\mbox{div}(A\nabla v) + B\nabla v=0,& \ \ \mbox{in $B^-_{\rho_2}(0)$},
\\
A(y',0)\nabla v\cdot\nu'=0,&\ \  \mbox{on $B'_{\rho_2}(0)$}.
\end{array}
\right.
\end{equation}

Hence we are under the assumptions of Theorem 1.3 in \cite{ae} and thus
we
can infer that there exists a radius $\rho_3, 0<\rho_3<\rho_2$,
depending on the \emph{a priori data} only, such that

\begin{eqnarray}\label{doublingperv}
\int_{B^-_{\beta \rho}(0)}|v|^2 \le c{\b}^K\int_{B^-_{\rho}(0)}|v|^2\ ,
 \end{eqnarray}
for every $\rho,\beta$ such that $\b>1$ and $0<\b \rho\le\rho_3$,
where $c>0$ is constant depending on the \emph{a priori data} only,
and $K>0$ depends on the \emph{a priori data} and increasingly on
 \begin{eqnarray}\label{Nvero}
N(\rho_3)=\rho_3\frac{\int_{B^-_{\rho_3}(0)}A\nabla
v\cdot\nabla \bar{v}+Re(\bar{v}\ \mbox{div}(A\nabla v))}{\int_{\partial
B^-_{\rho_3}(0)\setminus B'_{\rho_3(0)}}\mu |v|^2}\ ,
\end{eqnarray}
where we denote
\begin{eqnarray}
\mu(x)=\frac{A(x)x\cdot x}{|x|^2},\  \ \mbox{for every}\ x \in
B^-_{\rho_2}(0).
\end{eqnarray}
By \eqref{adoesc} it follows that
\begin{eqnarray}
c_3\le\mu(x)\le c_4,\ \ \mbox{for every}\ x \in B^-_{\rho_2}(0).
\end{eqnarray}
Let us observe that the proof of Theorem 1.3 in \cite{ae} needs, in
this
context, a slight modification due to the fact that we deal with
complex valued functions. We omit the details.

Denoting by
 \begin{eqnarray}\label{N}
\tilde{N}(\rho_3)=\frac{\int_{B^-_{\rho_3}(0)}\rho^2_3|\nabla
v|^2+|v|^2}{\int_{B^-_{\rho_3}(0)}|v|^2}\ ,
\end{eqnarray}
we notice, following the arguments in \cite[Lemma 3.3]{am}, that
 \begin{eqnarray}
N(\rho_3)\le C \tilde{N}(\rho_3),
\end{eqnarray}
where $C>0$ is a constant depending on the \emph{a priori data} only.

By \eqref{inclusioni2}, it follows, that for every $r$ and $\b>1$ such
that $0<r<\b r<\frac{\rho_3}{2}$
\begin{eqnarray}\label{w}
\int_{\G_{I,\b r}(0)}|z|^2
\le C\int_{B^-_{2\beta r}(0)}|v|^2\ ,
 \end{eqnarray}
where $C>0$ is a constant depending on $r_0, M, \Lambda$ only.
Moreover, by \eqref{doublingperv} and by \eqref{inclusioni2} we have
that
\begin{eqnarray}\label{ww}
\int_{B^-_{2\beta r}(0)}|v|^2 \le C(2\b
c_1)^K\int_{B^-_{\frac{r}{c_1}}(0)}|v|^2\le C(2\b c_1)^K\int_{\G_{I,
r}(0)}|z|^2,
\end{eqnarray}
where $C>0$ is a constant depending on $r_0, M, \Lambda$ only.

Combining \eqref{w} and \eqref{ww}, we have that
\begin{eqnarray}
\int_{\G_{I,\b r}}|z|^2\le C(2\b c_1)^K\int_{\G_{I, r}(0)}|z|^2\ .
\end{eqnarray}
Finally the last inequality, \eqref{uno},\eqref{psize} imply that
\begin{eqnarray}
\int_{\G_{I,\b r}}|u|^2\le C(\b)^K\int_{\G_{I, r}(0)}|u|^2\ ,
\end{eqnarray}
where $C>0, K>0$ are constants depending on \emph{a priori data} and
on $\tilde{N}(\rho_3)$ only. Thus the Lemma follows with
\begin{eqnarray}\label{raggio}
\bar{r}=\frac{\rho_3}{2}.
\end{eqnarray}
It only remains to majorize the quantity \eqref{N} by a constant
depending on the \emph{a priori data} only. Let us observe that by
\eqref{inclusioni2}, by \eqref{uno} and by \eqref{psize}, we have
that
\begin{eqnarray}
\int_{B^-_{\rho_3}(0)}|\nabla v|^2+|v|^2\le C
\int_{\G_{I,\rho_3c_1}(0)}|\nabla u|^2+|u|^2,
\end{eqnarray}
where $C>0$ is a constant depending on the \emph{a priori data} only.
Moreover, by the above inequality and by \eqref{ciunoalfa}, we can
conclude that
\begin{eqnarray}\label{nnn}
\int_{B^-_{\rho_3}(0)}|\nabla v|^2+|v|^2\le C,
\end{eqnarray}
where $C>0$ is a constant depending on \emph{a priori data} only.

On the other hand, we have that choosing
$P_0=\frac{M}{8\sqrt{1+M^2}\rho_3}\nu$ and
$\rho_4=\frac{1}{32}\frac{M}{\sqrt{1+M^2}}\rho_3$, where
$\nu$ is the outer unit normal to $D$ at $0$, it follows that
$B_{\rho_4}(P_0)\subset \G_{I,\frac{\rho_3}{2}}(0)$.

Thus, by \eqref{inclusioni2} and by \eqref{psize} it follows that
\begin{eqnarray}\label{nn}
\int_{B^-_{\rho_3}(0)}|v|^2\ge C
\int_{\G_{I,\frac{\rho_3}{2}}(0)}|u|^2\ge C
\int_{B_{\rho_4}(P_0)}|u|^2\ ,
\end{eqnarray}
where $C>0$ is a constant depending on the \emph{a priori data} only.

 Let us consider a point $Q\in
\R^3\setminus D^+_{2R_0}$ such that
\begin{eqnarray}\label{sub}
B_{4\rho_4}(Q)\subset \R^3\setminus \overline{D}^+_{2R_0},
\end{eqnarray}
where $R_0$ is the radius introduced in Corollary \ref{lb}.
Dealing as in the proof of Theorem \ref{stabpc}, we cover a path
joining $P_0$ to
$Q$ by a chain of balls of radius $\rho_4$ pairwise tangent to each
other. Hence, by an iterated use of the three spheres inequality, we
have that the following holds
\begin{eqnarray}
\|u\|_{L^2(B_{\frac{\rho_4}{4}(Q)})}\le C
\|u\|^{{\tau}^s}_{L^2(B_{\rho_4(P_0)})}\ ,
\end{eqnarray}
where $C>0, s>0$ and $\tau, 0<\tau<1$ are constants depending on the
\emph{a priori data} only.
By the last inequality, by \eqref{sub} and by \eqref{lowerbound}, we
can
infer that
\begin{eqnarray}\label{e}
\|u\|_{L^2(B_{\rho_4(P_0)})}\ge \left(\frac{\pi\rho^3_4}{C 48}
\right)^{\frac{1}{{\tau}^s}}.
\end{eqnarray}
Hence, by \eqref{e} and by \eqref{nn}, we have that
\begin{eqnarray}\label{nnnn}
\int_{B^-_{\rho_3}(0)}|v|^2\ge C,
\end{eqnarray}
where $C>0$ is a constant depending on \emph{a priori data} only.
Hence, by \eqref{nnn} and by \eqref{nnnn}, we can majorize
$\tilde{N}(\rho_3)$ by
a constant depending on the \emph{a priori data} only and thus the
Lemma follows.
\end{proof}

\begin{theorem}[Surface doubling inequality]\label{di}
Let $u$ be the solution to the problem \eqref{Sc}, then there exists a
constant $C>0$ depending on the \emph{a priori data} only such that,
for
every $x_0\in\G^{r_0}_{I}$ and for every $r\in (0,\frac{\bar{r}}{4})$,
the following holds
\begin{eqnarray}\label{doublingon}
\int_{\Delta_{I, 2r}(x_0)}|u|^2 \mbox{d}\sigma\le C\int_{\Delta_{I,
r}(x_0)}|u|^2\mbox{d}\sigma\ .
 \end{eqnarray}
\end{theorem}
\begin{proof}
Let $x_0\in\G^{r_0}_I$ and let $z\in H^1(\G_{I,r_1}(x_0))$ and
$\bar{r}$ be, respectively, the
solution to the problem \eqref{z2} defined by \eqref{defz} and the
radius introduced in \eqref{raggio}. By a regularity estimate at the
boundary, (see for instance \cite[p.777]{abrv}) we
have that, for any $r\in (0,\frac{\bar{r}}{4})$, the following holds
\begin{eqnarray}\label{int}
\int_{\Delta_{I,r}(x_0)}|\nabla_{t}z|^2\le
C\left(\frac{1}{r}\int_{\G_{I,2r}(x_0)}|\nabla
z|^2\right)^{1-\gamma}\left(\frac{1}{r^2}\int_{\Delta_{I,r}(x_0)}|z|^2\right)^{\gamma},
\end{eqnarray}
where $C>0$ and $0<\gamma<1$ are constants depending on the
\emph{a priori data} only and where $\nabla_{\mbox{t}}z$ represents the
tangential
gradient.

Thus, by the Young inequality we have that for every $\eps>0$ the
following holds
\begin{eqnarray}\label{schwartz}
\int_{\Delta_{I,r}(x_0)}|\nabla_{t}z|^2\le
\frac{C{\eps}^{\frac{1}{1-\gamma}}}{r}\int_{\G_{I,2r}(x_0)}|\nabla
z|^2+\frac{C}{{\eps}^{\frac{1}{\gamma}}r^2}
\int_{\Delta_{I,r}(x_0)}|z|^2,
\end{eqnarray}
where $C>0$ is a  constant depending on the
\emph{a priori data} only.

Moreover, by a well-known estimate of stability for the Cauchy problem
(see for instance \cite{tr}), we have
that
\begin{eqnarray}\label{grad}
\int_{\G_{I,\frac{r}{2}}(x_0)}|z|^2&&\le
Cr\left(\int_{\Delta_{I,r}(x_0)}|z|^2+r^2\int_{\Delta_{I,r}(x_0)}|\nabla_{t}z|^2\right)^{1-\delta}\cdot\\&&\cdot\left(\int_{\Delta_{I,r}(x_0)}|z|^2+r^2\int_{\Delta_{I,r}(x_0)}|\nabla_{t}z|^2+r\int_{\G_{I,r}(x_0)}|\nabla
z|^2\right)^{\delta},\nonumber
\end{eqnarray}
where $C>0$ and $0<\delta<1$ are constants depending on the
\emph{a priori data} only.

Hence, by \eqref{grad} and by the Young inequality, we have that for
every  $\beta>0$ the following holds
\begin{eqnarray}\label{gradsch}
\int_{\G_{I,\frac{r}{2}}(x_0)}|z|^2\le
\frac{C}{\eps^{\frac{\beta}{1-\delta}}}\left(r\int_{\Delta_{I,r}(x_0)}\!\!\!\!\!|z|^2+r^3\int_{\Delta_{I,r}(x_0)}|\nabla_{t}z|^2\right)+\\
+C\eps^{\frac{\beta}{\delta}}\left(r\int_{\Delta_{I,r}(x_0)}|z|^2+r^3\int_{\Delta_{I,r}(x_0)}|\nabla_{t}z|^2+r^2\int_{\G_{I,r}(x_0)}|\nabla
z|^2\right),\nonumber
\end{eqnarray}
where $C>0$ is a  constant depending on the \emph{a priori data}
only.

Choosing  $\beta$ in \eqref{gradsch} such that
$\beta=\frac{1-\delta}{1-\gamma}\gamma$ and inserting
\eqref{schwartz} in \eqref{gradsch}, we obtain
\begin{eqnarray}\label{gradsch2}
\int_{\G_{I,\frac{r}{2}}(x_0)}|z|^2\le
\frac{Cr}{{\eps}^{\frac{\gamma^2+1-\gamma}{\gamma(1-\gamma)}}}\int_{\Delta_{I,r}(x_0)}|z|^2
+C\eps r^2\int_{\G_{I,2r}(x_0)}|\nabla z|^2,\nonumber
\end{eqnarray}
where $C>0$ is a  constant depending on the
\emph{a priori data} only.
By the Caccioppoli inequality we have that
\begin{eqnarray}\label{gradsch3}
\int_{\G_{I,\frac{r}{2}}(x_0)}|z|^2\le
\frac{Cr}{{\eps}^{\frac{\gamma^2+1-\gamma}{\gamma(1-\gamma)}}}\int_{\Delta_{I,r}(x_0)}|z|^2
+C\eps \int_{\G_{I,4r}(x_0)}|z|^2,\nonumber
\end{eqnarray}
where $C>0$ is a  constant depending on the
\emph{a priori data} only.

Thus by \eqref{uno} and \eqref{psize} we can infer that
\begin{eqnarray}\label{gradsch4}
\int_{\G_{I,r}(x_0)}|u|^2\le
\frac{Cr}{{\eps}^{\frac{\gamma^2+1-\gamma}{\gamma(1-\gamma)}}}\int_{\Delta_{I,2r}(x_0)}|u|^2
+C\eps \int_{\G_{I,8r}(x_0)}|u|^2,\nonumber
\end{eqnarray}
where $C>0$ is a  constant depending on the
\emph{a priori data} only.

By \eqref{doublingat} it follows that
\begin{eqnarray}\label{sc3}
\int_{\G_{I,\frac{r}{2}}(x_0)}|u|^2\le
\frac{Cr}{{\eps}^{\frac{\gamma^2+1-\gamma}{\gamma(1-\gamma)}}}\int_{\Delta_{I,r}(x_0)}|u|^2
+C(8)^K\eps \int_{\G_{I,\frac{r}{2}}(x_0)}|u|^2,
\end{eqnarray}
where $C>0$ is a constant depending on the \emph{a priori data}
only.

Hence, choosing $\eps$ in \eqref{sc3} such that
$\eps=\frac{1}{2C(8)^K}$, we obtain that
\begin{eqnarray}\label{sc4}
\int_{\G_{I,\frac{r}{2}}(x_0)}|u|^2\le{Cr}\int_{\Delta_{I,r}(x_0)}|u|^2,
\end{eqnarray}
where $C>0$ is a constant depending on the \emph{a priori data}
only.

By applying again \eqref{doublingat} on the left hand side of
\eqref{sc4}, we obtain that
\begin{eqnarray}\label{sc5}
\int_{\G_{I,2r}(x_0)}|u|^2\le{Cr}\int_{\Delta_{I,r}(x_0)}|u|^2,
\end{eqnarray}
where $C>0$ is a constant depending on the \emph{a priori data}
only.

 Moreover, by a standard Dirichlet trace inequality, we have that
\begin{eqnarray}\label{sc6}
\int_{\Delta_{I,2r}(x_0)}|u|^2\le C\int_{\Delta_{I,r}(x_0)}|u|^2,
\end{eqnarray}
where $C>0$ is a constant depending on the \emph{a priori data}
only.

\end{proof}

\begin{corollary}[$A_{p}$ property on the boundary]\label{ap}
Let $u$ be the solution to the problem \eqref{Sc}, then there exist
$p>1, A>0$ constants depending on
the \emph{a priori data} only, such that, for every
$x_0\in\G^{r_0}_{I}$
and every $r\in (0,\frac{\bar{r}}{4})$, the following holds
\begin{eqnarray}\label{approp}
\!\!\!\!\!\!\!\!\left(\frac{1}{|\Delta_{I,r}(x_0)|}\int_{\Delta_{I,r}(x_0)}|u|^2
\mbox{d}\sigma
\right)\left(\frac{1}{|\Delta_{I,r}(x_0)|}\int_{\Delta_{I,r}(x_0)}|u|^{-\frac{2}{p-1}}
\mbox{d}\sigma \right)^{p-1}\le A .
\end{eqnarray}
\end{corollary}
\begin{proof}
Let $x_0\in\G^{r_0}_{I}$ and let $r\in (0,\frac{\bar{r}}{4})$, then by
a trace inequality, (see for instance \cite{Adams}, Chap. 5), it
follows that
\begin{eqnarray}
\|u\|_{L^4(\Delta_{I,r}(x_0))}\le C \|u\|_{H^1(\G_{I,r}(x_0))},
\end{eqnarray}
where $C>0$ is a constant depending on the \emph{a priori data} only.
By the Caccioppoli inequality we deduce that
\begin{eqnarray}\label{apcac}
\|u\|_{L^4(\Delta_{I,r}(x_0))}\le
\frac{C}{r}\|u\|_{L^2(\G_{I,2r}(x_0))}.
\end{eqnarray}
Applying the Doubling inequality \eqref{doublingat} on the right hand
side of \eqref{apcac}, we obtain that
\begin{eqnarray}\label{2}
\|u\|_{L^4(\Delta_{I,r}(x_0))}\le
\frac{C}{r}\|u\|_{L^2(\G_{I,r}(x_0))},
\end{eqnarray}
where $C>0$ is a constant depending on the \emph{a priori data} only.
Combining \eqref{sc4} and \eqref{2} we have that
\begin{eqnarray}
\|u\|_{L^4(\Delta_{I,r}(x_0))}\le
\frac{C}{\sqrt{r}}\|u\|_{L^2(\Delta_{I,
2r}(x_0))},
\end{eqnarray}
where $C>0$ is a constant depending on the \emph{a priori data} only.
Thus by the doubling inequality \eqref{doublingon} we have
\begin{eqnarray}
\|u\|_{L^4(\Delta_{I,r}(x_0))}\le
\frac{C}{\sqrt{r}}\|u\|_{L^2(\Delta_{I,
r}(x_0))}.
\end{eqnarray}
Hence, we infer that for every $r\in(0,\frac{\bar{r}}{4})$ and for
every $x_0\in\G^{r_0}_I$, the following holds
$$\left( \frac{1}{r^2}\int_{\Delta_{I,r}}|u|^4
\right)^{\frac{1}{4}}\le\left(
\frac{C}{r^2}\int_{\Delta_{I,r}}|u|^2
\right)^{\frac{1}{2}}, $$
obtaining a reverse H\"older inequality.

The result in \cite{cf} assures the existence of some $p>1$ and $A>0$
depending on the \emph{a priori data} only such that \eqref{approp}
holds.
\end{proof}

\begin{proofl}
Let $x_0$ be a point in $\G^{r_0}_{I}$. Let us pick
$r=\frac{\bar{r}}{8}$, thus by \eqref{sc4} with $u=u_2$ it follows that
\begin{eqnarray}\label{ll}
\int_{\Delta_{I,\frac{\bar{r}}{8}}(x_0)}|u_2|^2
\mbox{d}\sigma\ge
C\int_{\G_{I,\frac{\bar{r}}{16}(x_0)}}|u_2|^2\mbox{d}x,
\end{eqnarray}
where $C>0$ is a constant depending on the \emph{a priori data} only.

Let $P_0$ and $\rho_4>0$ be, respectively a point and a radius, such
that
$B_{\rho_4}(P_0)\subset\G_{I,\frac{\bar{r}}{16}(x_0)}$. By rephrasing
the argument leading to \eqref{e} we deduce by \eqref{ll} that
\begin{eqnarray}\label{lll}
\int_{\Delta_{I,\frac{\bar{r}}{8} }(x_0)}|u_2|^2 \mbox{d}\sigma\ge C,
\end{eqnarray}
where $C>0$ is a constant depending on the \emph{a priori data} only.

Combining \eqref{approp} and \eqref{lll}, we have that for every $x_0
\in \G^{r_0}_{I}$ the following holds
\begin{eqnarray}\label{es}
\left(\int_{\Delta_{I,\frac{\bar{r}}{8}}(x_0)}|u_2|^{-\frac{2}{p-1}}
\mbox{d}\sigma \right)^{p-1}\le C,
\end{eqnarray}
where $C>0$ is a constant depending on the \emph{a priori data} only.

Let us now consider $x\in \Delta_{I,\frac{\bar{r}}{8}}(x_0)$, then it
follows that
\begin{eqnarray*}
|\lambda_1(x)
-\lambda_2(x)|&=&\left|-\lambda_1(x)\frac{u_1(x)-u_2(x)}{u_2(x)}+\frac{1}{i
u_2(x)}
\left(\displaystyle{\frac{\partial u_2(x)}{\partial
\nu}}-\displaystyle{\frac{\partial u_1(x)}{\partial \nu}} \right)
 \right|\le \\
&\le&|\lambda_1(x)|\frac{|u_1(x)-u_2(x)|}{|u_2(x)|} +\frac{1}{|
u_2(x)|}\left|\displaystyle{\frac{\partial u_2(x)}{\partial
\nu}}-\displaystyle{\frac{\partial u_1(x)}{\partial \nu}} \right|_{\ .}
\end{eqnarray*}
Then by Theorem \ref{stabpc} and by \eqref{L} we have that, if
$0<\eps<\eps_0$, then
\begin{eqnarray}\label{sp}
|\lambda_1(x)
-\lambda_2(x)|\le (\Lambda +1)\eta(\eps)\frac{1}{|
u_2(x)|}.
\end{eqnarray}
Hence denoting by $\delta=\frac{2}{p-1}$, \eqref{sp} yields to
\begin{eqnarray}\label{sp2}
\left(\int_{\Delta_{I,\frac{\bar{r}}{8}(x_0)}}|\lambda_1(x)
-\lambda_2(x)|^{\delta}\right)^{\frac{1}{\delta}}\le (\Lambda
+1)\eta(\eps)\left(\int_{\Delta_{I,\frac{\bar{r}}{8}(x_0)}}\frac{1}{|
u_2(x)|^{\delta}}\right)^{\frac{1}{\delta}}.
\end{eqnarray}
By \eqref{es} and by a possible replacement of the constant $C$ in
\eqref{eta}, we have that
\begin{eqnarray}\label{sp3}
\left(\int_{\Delta_{I,\frac{\bar{r}}{8}(x_0)}}|\lambda_1(x)
-\lambda_2(x)|^{\delta}\right)^{\frac{1}{\delta}}\le\eta(\eps).
\end{eqnarray}
By the a priori bound \eqref{L}, we can infer that
\begin{eqnarray}\label{sp4}
|\lambda_1(x)-\lambda_2(x)|\le
|\lambda_1(x)-\lambda_2(x)|^{\frac{\delta}{2}}(2\Lambda)^{1-\frac{\delta}{2}}\
.
\end{eqnarray}
Integrating the above inequality with respect to $x$ over
$\Delta_{I,\frac{\bar{r}}{8}(x_0)}$ we have
\begin{eqnarray}\label{sp5}
\|\lambda_1(x)-\lambda_2(x)\|_{L^2(\Delta_{I,\frac{\bar{r}}{8}(x_0)})}\le
(2\Lambda)^{1-\frac{\delta}{2}}\left(\int_{\Delta_{I,\frac{\bar{r}}{8}(x_0)}}|\lambda_1(x)
-\lambda_2(x)|^{\delta} \right)^{\frac{1}{2}}.
\end{eqnarray}

Hence, by  a possible further replacement of the constants $C,\theta$
in \eqref{eta}, we can infer that the last inequality and \eqref{sp3} yield to
\begin{eqnarray}\label{sp6}
\|\lambda_1(x)-\lambda_2(x)\|_{L^2(\Delta_{I,\frac{\bar{r}}{8}(x_0)})}\le
\eta(\eps)\ .
\end{eqnarray}
By an interpolation inequality, see for instance \cite[p.777]{abrv},
we have that
\begin{eqnarray}\label{sp7}
\|\lambda_1-\lambda_2\|_{L^{\infty}(\Delta_{I,\frac{\bar{r}}{8}(x_0)})}\le
C\|\lambda_1-\lambda_2\|^{\frac{1}{2}}_{L^{2}(\Delta_{I,\frac{\bar{r}}{8}(x_0)})}\|\lambda_1-\lambda_2\|^{\frac{1}{2}}_{C^{0,1}(\Delta_{I,\frac{\bar{r}}{8}(x_0)})},
\end{eqnarray}
where $C>0$ is a constant depending on the \emph{a priori data} only.
Hence by \eqref{L}, it follows that
\begin{eqnarray}\label{sp8}
\|\lambda_1-\lambda_2\|_{L^{\infty}(\Delta_{I,\frac{\bar{r}}{8}(x_0)})}\le C(2\Lambda)^{\frac{1}{2}}
\|\lambda_1-\lambda_2\|^{\frac{1}{2}}_{L^{2}(\Delta_{I,\frac{\bar{r}}{8}(x_0)})}.
\end{eqnarray}

Combining \eqref{sp6} with \eqref{sp8} we obtain, by  a possible
further replacement of the constants $C,\theta$ in \eqref{eta}, that
\begin{eqnarray}\label{sp9}
\|\lambda_1-\lambda_2\|_{L^{\infty}(\Delta_{I,\frac{\bar{r}}{8}(x_0)})}\le\eta(\eps).
\end{eqnarray}

Let us cover $\G^{r_0}_I$ with the sets
$\Delta_{I,\frac{\bar{r}}{8}}(x_j),\ j=1,\dots,J$, with
$x_j\in\G^{r_0}_I$.

Let $i$ be an index such that
\begin{eqnarray}\label{sp10}
\|\lambda_1-\lambda_2\|_{L^{\infty}(\Delta_{I,\frac{\bar{r}}{8}(x_i)})}=\|\lambda_1-\lambda_2\|_{L^{\infty}(\G^{r_0}_{I})}.
\end{eqnarray}

Thus, by a further possible replacement of the constant
$C,\theta$ in \eqref{eta}, we deduce \eqref{aim} by combining
\eqref{sp10} and \eqref{sp9} with $x_0=x_i$.
\end{proofl}

\end{document}